\documentclass[a4paper,fleqn]{cas-sc}

\usepackage[authoryear,longnamesfirst]{natbib}
\usepackage{subcaption} 
\usepackage{siunitx}
\usepackage{ulem}

\def\tsc#1{\csdef{#1}{\textsc{\lowercase{#1}}\xspace}}
\tsc{WGM}
\tsc{QE}

\newcommand{\kmax}{K_{\max}} 
\newcommand{\lmax}{L_{\max}} 
\newcommand{\imax}{I_{\max}} 
\newcommand{\state}{\alpha} 
\newcommand{\basis}{\phi} 
\newcommand{\radbasis}{R} 
\newcommand{\angbasis}{Y} 
\newcommand{\intbasis}{H} 
\newcommand{\maxwellianref}{\tilde{M}} 
\newcommand{\intweightref}{\tilde{E}} 

\newcommand{\eq}{\mathrm{eq}} 
\newcommand{\densityref}{\tilde{\rho}} 
\newcommand{\tempref}{\tilde{T}} 

\newcommand{\vpre}{\mathbf{v}} 
\newcommand{\wpre}{\mathbf{w}} 
\newcommand{\vpost}{\mathbf{v}'} 
\newcommand{\wpost}{\mathbf{w}'} 

\newcommand{\vrelpre}{\mathbf{u}} 
\newcommand{\vrelpost}{\mathbf{u}'} 

\newcommand{\vpreunit}{\hat{\mathbf{v}}}    
\newcommand{\wpreunit}{\hat{\mathbf{w}}}    
\newcommand{\vpostunit}{\hat{\mathbf{v}}'}  
\newcommand{\scatterdir}{\hat{\mathbf{u}}'} 
\newcommand{\vreldir}{\hat{\mathbf{u}}} 

\newcommand{\Ipre}{I} 
\newcommand{\Jpre}{J} 
\newcommand{\Ipost}{I'} 
\newcommand{\Jpost}{J'} 
\newcommand{\Etot}{E} 

\newcommand{\polyalpha}{\nu} 
\newcommand{\dofpoly}{\delta} 

\newcommand{\partkin}{R} 
\newcommand{\partint}{r} 
\newcommand{\phipart}{\varphi_{\polyalpha}} 
\newcommand{\psipart}{\psi_{\polyalpha}} 

\newcommand{\incangle}{\beta}        
\newcommand{\defangle}{\chi}        
\newcommand{\aziangle}{\epsilon}    

\newcommand{\vmag}{v}               
\newcommand{\wmag}{w}               
\newcommand{\urel}{u}               

\newcommand{\Ctens}{C} 
\newcommand{\Rtens}{\mathcal{R}} 
\newcommand{\Gtens}{\mathcal{G}} 
\newcommand{\wignerD}{\mathcal{D}} 

\newcommand{\cmvar}{t} 
\newcommand{\srel}{s} 
\newcommand{\fracrel}{x} 
\newcommand{\cmenergy}{z} 
\newcommand{\halfangle}{y} 
\newcommand{\duffy}{h} 
\newcommand{\vhs}{\gamma} 
\newcommand{\chan}{c} 

\newcommand{\dof}{N} 
\newcommand{\numQ}{N_Q} 
\newcommand{\numT}{N_T} 
\newcommand{\numG}{N_G} 
\newcommand{\numS}{N_S} 
\newcommand{\inelprob}{\omega} 
\newcommand{\Kmono}{K_{\dofpoly}} 

\newcommand{\modexp}{\vartheta} 
\newcommand{\modexpf}{\vartheta_{\!f}} 
\newcommand{\modamp}{\hat{\eta}} 
\newcommand{\modampf}{\hat{\eta}_{f}} 
\newcommand{\auxvar}{\lambda} 
\newcommand{\auxmap}{\xi} 
\newcommand{\auxscale}{\sigma} 
\newcommand{\numaux}{N_{\auxvar}} 

\newcommand{\prandtl}{\mathrm{Pr}} 
\newcommand{\shearvisc}{\mu} 
\newcommand{\bulkvisc}{\mu_{b}} 
\newcommand{\heatcond}{\kappa} 

\newcommand{\vcm}{\bar{\mathbf{v}}} 
\newcommand{\polyblock}{\Omega} 
\newcommand{\stageM}{\mathcal{M}} 
\newcommand{\stageS}{\mathcal{S}} 
\newcommand{\stageN}{\mathcal{N}} 
\newcommand{\stageO}{\mathcal{O}} 
\newcommand{\radweight}{\mathcal{K}} 
\newcommand{\redker}{\tilde{B}} 

\usepackage{xcolor}

\begin{document}
\let\WriteBookmarks\relax
\def\floatpagepagefraction{1}
\def\textpagefraction{.001}

\shorttitle{Wigner-Eckart Factorization of the Polyatomic Boltzmann Collision Operator}    

\shortauthors{T. Ke\ss ler, B.W.T. Gieling, R.R. Hiemstra, and M.R.A. Abdelmalik}  

\title [mode = title]{Wigner-Eckart Factorization of the Polyatomic Boltzmann Collision Operator}  



\author[1,2]{Torsten Ke\ss ler}

\cormark[1]

\ead{t.kessler@tue.nl}

\credit{Conceptualization, Methodology, Software, Validation, Writing, Review \& Editing}

\author[1]{B.W.T. Gieling}

\ead{b.w.t.gieling@tue.nl}

\credit{Software, Validation, Review \& Editing}

\author[1,2]{René R. Hiemstra}[orcid=0000-0002-0179-9606]

\fnmark[1]

\ead{r.r.hiemstra@tue.nl}

\credit{Conceptualization, Methodology, Software, Validation, Writing, Review \& Editing}

\author[1,2]{Michael R.A. Abdelmalik}

\ead{m.r.a.abdelmalik@tue.nl}

\credit{Conceptualization, Methodology, Review \& Editing}

\affiliation[1]{organization={Department of Mechanical Engineering, Eindhoven University of Technology (TU/e)},
            country={The Netherlands}}
            
\affiliation[2]{organization={Simkinetic},
            url={www.simkinetic.tech},
            city={Veldhoven},
            country={The Netherlands}}

\cortext[1]{Corresponding author}

\fntext[1]{Marie Sk\l odowska-Curie Fellow.}


\begin{abstract}
We extend the Wigner--Eckart factorization of the spectral Boltzmann collision operator to polyatomic gases with continuous internal energy. Because internal energies are invariant under spatial rotations, the $\mathrm{SO}(3)$ reduction survives the Borgnakke--Larsen energy exchange, and the twelve-dimensional collision integral collapses onto a nine-dimensional kinematic core. The core splits into a sparse geometric tensor, evaluated exactly, and a dense physical tensor, integrated by singularity-resolving Gauss rules with an auxiliary Laplace representation of the fractional energy couplings. The quadrature attains near machine precision at the fractional exponents of real gases. The collision invariants are embedded exactly, preserving the translational--internal energy exchange. The factorization compresses the operator by three to nearly four orders of magnitude and accelerates its evaluation $40$-fold over dense formulations. The method is validated against the exact monatomic limit, Landau--Teller relaxation, and an analytic frozen-channel Prandtl number, and it matches a published calibration of the same kernel for N$_2$, CO, and H$_2$.
\end{abstract}


\begin{highlights}
\item Wigner-Eckart factorization extended to the polyatomic Borgnakke-Larsen model.
\item Collision integral reduced to a 9D kinematic core via SO(3) rotational symmetry.
\item Macroscopic angular geometry decouples from continuous internal energy exchange.
\item Tailored coordinate transformations enable spectral quadrature for the 9D operator.
\item Null-space embedding natively conserves mass, momentum, and total molecular energy.
\end{highlights}


\maketitle

\section{Introduction}

Rarefied molecular flows, such as hypersonic reentry aerothermodynamics \citep{ivanov:1998} and vacuum and microscale process flows \citep{sharipov:1998}, lie beyond the reach of continuum models and require the Boltzmann equation for polyatomic gases. The exchange between translational and internal energy that distinguishes these gases is macroscopically visible: bulk viscosities of common polyatomic gases exceed their shear viscosities by up to several orders of magnitude \citep{cramer:2012}. A continuum bulk-viscosity coefficient is an inadequate surrogate for finite-rate exchange in shock structure and acoustic absorption \citep{graves:1999}. This regime is dominated by the direct simulation Monte Carlo (DSMC) method \citep{bird:1994,boyd:2017}, which treats internal energy through phenomenological exchange models calibrated against measured relaxation data \citep{millikan:1963}. The statistical noise and slow convergence of this stochastic approach motivate deterministic alternatives.

The deterministic route is well developed for monatomic gases. Fourier-spectral methods offer fast, $N \log N$ evaluation of the collision operator \citep{pareschi:2000,mouhot:2006,gamba:2017} but sacrifice exact conservation, restored only by a posteriori correction \citep{gamba:2009}. Global spectral bases in the Burnett--Sonine lineage \citep{kumar:1966,cai:2020} conserve naturally and converge spectrally, but their dense Cartesian collision tensors incur $\mathcal{O}(L_{\max}^6)$ storage and arithmetic. These costs are addressed by conservative test-space design \citep{kessler:2019} and, most recently, by representation-theoretic and low-rank compression \citep{hanke:2023,hu:2022}. For polyatomic gases, deterministic methods are less developed: the collision integral grows from eight to twelve dimensions, and the internal energy couples into the kernel, the basis, the quadrature, and the conservation constraints at once. The closest prior solver treats discrete internal levels as separate species with optimization-based conservation \citep{munafo:2014}. It does not discretize the continuous-internal-energy operator. Model-equation surrogates and moment hierarchies \citep{andries:2000,rahimi:2016} avoid the collision integral altogether. A spectrally convergent, exactly conservative Galerkin discretization of the full continuous-energy polyatomic operator has, to our knowledge, not been available.

\subsection{Historical Context: From Molecular Gases to Computable Kinetics}
The Boltzmann equation describes structureless particles, but the gases of engineering interest are molecular: internal degrees of freedom store energy, exchange it in collisions, and alter transport. The classical account of these effects runs from Eucken's correction of the monatomic thermal conductivity \citep{eucken:1913} through the semi-classical equation of Wang Chang and Uhlenbeck \citep{wangchang:1951} and the transport theory of \citet{mason:1962}, and is codified in the two-volume treatise of \citet{mccourt:1990,mccourt:1991}. A computation-oriented line opened when \citet{borgnakke:1975} replaced the quantum state sum by a phenomenological redistribution of a continuous internal energy, the workhorse model of the direct simulation Monte Carlo (DSMC) method \citep{bird:1994}. \citet{bourgat:1994} showed that this picture supports a Boltzmann equation with microreversible collisions and an H-theorem, extended to polytropic gases by \citet{desvillettes:2005}; on this foundation rest Cauchy theory \citep{gamba2023cauchy,alonso:2025}, a unifying internal-state framework \citep{borsoni:2022}, and explicit collision kernels calibrated to experiment \citep{djordjic:2021pre,djordjic:2023}. In parallel, deterministic spectral methods for the collision operator have matured since the mid-1990s \citep{pareschi:2000}, with fast algorithms \citep{mouhot:2006,gamba:2017} and conservative, structure-exploiting formulations \citep{gamba:2009,kessler:2019,hanke:2023}, but almost exclusively for monatomic gases. The present work joins the two lines: a spectral-Galerkin discretization of the continuous-internal-energy collision operator with the extended DSMC kernel of \citet{djordjic:2023}, factorized through its $\mathrm{SO}(3)$ structure by the Wigner--Eckart theorem as in the monatomic companion \citep{hiemstra2026wigner}, exactly conservative, and calibrated against the measured transport coefficients of N$_2$, CO, and H$_2$.

\subsection{Core Contributions}
This paper extends the Wigner--Eckart factorization of the spectral Boltzmann collision operator \citep{hiemstra2026wigner} from monatomic to polyatomic gases. Specifically, it introduces the following primary contributions:

\begin{itemize}
  \item \textbf{Dimensional Reduction Extended to Internal States:} We carry the Wigner--Eckart factorization from the eight-dimensional monatomic collision integral to the twelve-dimensional integral of the Borgnakke--Larsen/Bourgat model. Because the internal energies and partition parameters are rotational scalars, the $\mathrm{SO}(3)$ structure of the collision geometry survives intact: the operator separates into the identical sparse geometric Gaunt tensor and a rotationally invariant physical tensor over a nine-dimensional kinematic core.
  \item \textbf{Spectrally Convergent Quadrature for the Extended DSMC Kernel:} We compose the quadrature so that every integration direction is exact or geometrically convergent by construction, introducing an auxiliary Laplace representation of the fractional total-energy couplings on both collision channels. The integer-exponent limits reach machine precision. At the fractional exponents of the production gases a slow algebraic endpoint tail remains, far below working accuracy.
  \item \textbf{Exact Conservation by Null-Space Embedding:} We extend the monatomic null-space embedding to the polyatomic invariants, where the total energy spans a two-mode pair rather than a single test row. An orthogonal projection removes only the invariant combination and leaves the physical translational--internal (Landau--Teller) exchange mode exactly as computed. Mass is conserved exactly and total energy to machine precision under nonlinear evolution.
  \item \textbf{Compressed Storage and Contraction on the Enlarged Radial--Internal Block:} The polyatomic construction enlarges the dense block of the operator along the internal axis, and the monatomic storage and contraction strategies apply verbatim to it. The compression over the dense Cartesian tensor grows along both truncation axes, and the cache-optimized angular-first contraction evaluates the operator $40$-fold faster than the dense baseline.
  \item \textbf{Validation Against Analytic Limits and an Independent Calibration:} The discrete operator reproduces an independent monatomic build at the internal collapse point, the Wang Chang--Uhlenbeck eigenvalue spectrum, the analytic Landau--Teller relaxation rates, the closed-form frozen Prandtl number of \citet{djordjic:2023}, and the classical Chapman--Enskog hard-sphere viscosity factor. These comparisons are parameter-free. At the published kernel parameters it reproduces the measured transport coefficients of N$_2$, CO, and H$_2$ with nothing fitted, but that agreement with experiment is inherited from the calibration of those parameters rather than predicted. What the comparison establishes is a code-to-code verification: our nine-dimensional quadrature reproduces their symbolically evaluated production coefficients to within $0.0002\%$.
\end{itemize}

\subsection{Outline of the Paper}
Section~\ref{sec:background} recalls the continuous-internal-energy collision model, its micro-reversible weak form, and the tensorized spectral Galerkin basis. Section~\ref{sec:reduction} derives the reduction of the twelve-dimensional collision integral to a nine-dimensional kinematic core and the Wigner--Eckart factorization. Section~\ref{sec:numerical} details the singular quadrature architecture, the auxiliary Laplace treatment of the extended kernel, the conservation embedding, and the tensor contraction strategies. Section~\ref{sec:results} validates the method: quadrature convergence, the monatomic limit, conservation, transient relaxation, transport coefficients for N$_2$, CO, and H$_2$ against experiment, and computational performance. Conclusions are drawn in Section~6.

\section{Background and Notation}
\label{sec:background}

This section defines the continuous Borgnakke--Larsen collision operator, its invariant measure and asymmetric weak form, the tensorized spectral Galerkin basis, and the discrete collision tensor.

Throughout this paper, boldface letters denote vectors in $\mathbb{R}^3$ (e.g., $\vpre, \wpre$), their non-bold counterparts denote scalar magnitudes or speeds (e.g., $\vmag = |\vpre|, \wmag = |\wpre|$), and hats denote unit directional vectors (e.g., $\vpreunit = \vpre/\vmag$). Primed variables designate post-collision states (e.g., $\vpost, \Ipost$).

\subsection{Continuous Polyatomic Collision Model}

We consider the initial value problem for the spatially homogeneous Boltzmann equation,
\begin{equation}
    \frac{\partial f(t,\vpre,\Ipre)}{\partial t} = Q(f,f)(t,\vpre,\Ipre), \quad t \in \mathbb{R}^+, \quad \vpre \in \mathbb{R}^3, \quad \Ipre \in [0, \infty),
\end{equation}
for the distribution function $f$. Each particle state carries a translational velocity $\vpre$ and a microscopic internal energy $\Ipre$.

For a colliding pair with pre-collision states $(\vpre, \Ipre)$ and $(\wpre, \Jpre)$, the relative velocity is $\vrelpre = \vpre - \wpre$ with magnitude $\urel$. Collisions conserve total momentum and total molecular energy. The total collision energy is
\begin{equation}
    \Etot = \frac{1}{4}\urel^2 + \Ipre + \Jpre = \frac{1}{4}(\urel')^2 + \Ipost + \Jpost.
\end{equation}

Following the Borgnakke--Larsen parameterization \citep{borgnakke:1975}, the collision is parameterized by the scattering direction $\scatterdir \in \mathbb{S}^2$ and two scalar partition parameters $\partkin, \partint \in [0,1]$. The parameter $\partkin$ dictates the fraction of the total energy converted into kinetic energy, yielding the post-collision relative velocity $\vrelpost = 2\sqrt{\partkin \Etot} \scatterdir$, which establishes the polar deflection angle $\cos \defangle = \vreldir \cdot \scatterdir$. The parameter $\partint$ distributes the remaining internal energy between the two interacting particles. The post-collision states are uniquely determined by
\begin{subequations}
\begin{align}
    \vpost &= \frac{1}{2}(\vpre + \wpre + \vrelpost), \quad \wpost = \frac{1}{2}(\vpre + \wpre - \vrelpost), \\
    \Ipost &= \partint(1-\partkin)\Etot, \quad \Jpost = (1-\partint)(1-\partkin)\Etot.
\end{align}
\end{subequations}

In the continuous-internal-energy framework of \citet{bourgat:1994}, with the polytropic weight of \citet{desvillettes:2005}, the polyatomic collision operator is the multi-dimensional integral
\begin{equation}
    Q(f,f)(\vpre,\Ipre) = \int_{\mathbb{R}^3} \int_0^\infty \int_0^1 \int_0^1 \int_{\mathbb{S}^2} B \left[ f(\vpost,\Ipost)f(\wpost,\Jpost) \left(\frac{\Ipre\Jpre}{\Ipost\Jpost}\right)^\polyalpha - f(\vpre,\Ipre)f(\wpre,\Jpre) \right] d\Phi,
\end{equation}
where $B \equiv B(\urel, \cos\defangle, \Ipre, \Jpre, \partkin, \partint)$ is the physical transition probability rate. The differential measure $d\Phi = \phipart(\partint) \psipart(\partkin) (1-\partkin) \partkin^{1/2} \, d\scatterdir \, d\partkin \, d\partint \, d\Jpre \, d\wpre$ incorporates the Jacobian of the collision transformation. The partition functions $\phipart(\partint) = (\partint(1-\partint))^\polyalpha$ and $\psipart(\partkin) = (1-\partkin)^{2\polyalpha}$ dictate the energy exchange. The constant $\polyalpha = \dofpoly/2 - 1$, with $\dofpoly$ the number of internal degrees of freedom of the gas, follows the convention of \citet{djordjic:2023}.

\subsection{Micro-Reversibility and the Asymmetric Weak Form}

To ensure micro-reversibility, the collision integral relies on the invariant measure $dA$ \citep{gamba2023cauchy}, which absorbs the pre-collision density of states $(\Ipre\Jpre)^\polyalpha$ to neutralize the asymmetry of the collision transformation:
\begin{equation}
    dA = B \cdot \phipart(\partint) \psipart(\partkin) (1-\partkin) \partkin^{1/2} (\Ipre\Jpre)^\polyalpha \, d\partkin \, d\partint \, d\scatterdir \, d\Jpre \, d\wpre \, d\Ipre \, d\vpre.
\end{equation}

Multiplying the strong form by a test function $\zeta(\vpre, \Ipre)$ and integrating over the target state space yields the asymmetric weak form of the collision operator:
\begin{equation}
    \int_{\mathbb{R}^3} \int_{\mathbb{R}^+} Q(f,f) \zeta \, d\Ipre \, d\vpre = \int_{(\mathbb{R}^3)^2 \times (\mathbb{R}^+)^2 \times [0,1]^2 \times \mathbb{S}^2} \frac{f(\vpre,\Ipre)f(\wpre,\Jpre)}{(\Ipre\Jpre)^\polyalpha} \left[ \zeta(\vpost,\Ipost) - \zeta(\vpre,\Ipre) \right] dA.
\end{equation}
Substituting $dA$ cancels the asymmetric quotient $(\Ipre\Jpre)^\polyalpha$ analytically, reducing the integral to the plain transition measure:
\begin{equation}
    \int_{\mathbb{R}^3} \int_{\mathbb{R}^+} Q(f,f) \zeta \, d\Ipre \, d\vpre = \int \Big[ \zeta(\vpost,\Ipost) - \zeta(\vpre,\Ipre) \Big] f(\vpre,\Ipre)f(\wpre,\Jpre) \, B \cdot \phipart(\partint) \psipart(\partkin) (1-\partkin) \partkin^{1/2} \, d\partkin \dots d\vpre.
    \label{eq:weak_form_cancelled}
\end{equation}

\subsection{Tensorized Spectral Galerkin Bases}
\label{sec:basis}

We tensorize a three-dimensional velocity-space basis (associated Laguerre polynomials and spherical harmonics) with orthogonal polynomials for the internal energy space. With the composite spectral state index $\state = (k, l, m, i)$, the basis functions decouple the radial, angular, and internal energy dependencies:
\begin{equation}
    \basis_\state(\vpre, \Ipre) = \radbasis_{k,l}(\vmag) \angbasis_{l,m}(\vpreunit) \intbasis_i(\Ipre).
\end{equation}

The distribution function is expanded using time-dependent spectral coefficients $c_\state(t)$ weighted by the absolute equilibrium states:
\begin{equation}
    f(t,\vpre,\Ipre) \approx \maxwellianref(\vmag) \intweightref(\Ipre) \sum_{\state=1}^\dof c_\state(t) \basis_\state(\vpre,\Ipre),
\end{equation}
where $\dof$ is the total number of basis functions. The reference Maxwellian weight $\maxwellianref(\vmag)$ is defined by the reference density $\densityref$ and temperature $\tempref$:
\begin{equation}
    \maxwellianref(\vmag) = \frac{\densityref}{(2\pi \tempref)^{3/2}} \exp\left(-\frac{\vmag^2}{2\tempref}\right).
\end{equation}
The equilibrium internal-energy weight $\intweightref(\Ipre)$ contains the $\Ipre^\polyalpha$ density of states:
\begin{equation}
    \intweightref(\Ipre) = \frac{1}{\tempref^{\polyalpha+1} \Gamma(\polyalpha+1)} \Ipre^\polyalpha \exp\left(-\frac{\Ipre}{\tempref}\right).
\end{equation}

We nondimensionalize with the reference temperature $\tempref$, measuring velocities in units of the reference thermal speed $\sqrt{2\tempref}$, so that the translational energy $\vmag^2$ and the internal energies $\Ipre$ and $\Jpre$ are in units of $\tempref$. The same symbols denote the scaled variables from here on. In these units the two reference weights take their canonical orthogonality forms,
\begin{equation}
    \maxwellianref(\vmag) = \frac{\densityref}{\pi^{3/2}} e^{-\vmag^2}, \qquad
    \intweightref(\Ipre) = \frac{1}{\Gamma(\polyalpha+1)} \Ipre^\polyalpha e^{-\Ipre},
    \label{eq:reference_units}
\end{equation}
against which the radial polynomials $\radbasis_{k,l}$ and the internal polynomials $\intbasis_i$ are orthonormal, and which the quadrature rules of Section~\ref{sec:numerical} carry as their native weights.

\subsection{Macroscopic Collision Tensor}

Setting the test function $\zeta = \basis_{\state_1}(\vpre, \Ipre)$ and substituting the spectral ansatz into the reduced weak form \eqref{eq:weak_form_cancelled} yields the coupled system of ordinary differential equations
\begin{equation}
    \frac{\partial c_{\state_1}}{\partial t} = \sum_{\state_2=1}^\dof \sum_{\state_3=1}^\dof \Ctens_{\state_1, \state_2, \state_3} c_{\state_2} c_{\state_3}.
\end{equation}

The macroscopic Cartesian collision tensor $\Ctens$ separates into a positive gain and negative loss component, $\Ctens = \Ctens^+ - \Ctens^-$. Because the reference weights $\intweightref(\Ipre)\intweightref(\Jpre)$ provide the $\Ipre^\polyalpha \Jpre^\polyalpha$ terms required for detailed balance, the discrete tensors assemble without micro-reversibility quotients.

The exact Cartesian tensors evaluate to twelve-dimensional integrals over the absolute laboratory frame:
\begin{subequations}
\begin{align}
    \Ctens_{\state_1, \state_2, \state_3}^+ &= \int_{(\mathbb{R}^3)^2 \times (\mathbb{R}^+)^2 \times [0,1]^2 \times \mathbb{S}^2} B \cdot \basis_{\state_1}(\vpost, \Ipost) \basis_{\state_2}(\vpre, \Ipre) \basis_{\state_3}(\wpre, \Jpre) \, d\tilde{\Phi}, \\
    \Ctens_{\state_1, \state_2, \state_3}^- &= \int_{(\mathbb{R}^3)^2 \times (\mathbb{R}^+)^2 \times [0,1]^2 \times \mathbb{S}^2} B \cdot \basis_{\state_1}(\vpre, \Ipre) \basis_{\state_2}(\vpre, \Ipre) \basis_{\state_3}(\wpre, \Jpre) \, d\tilde{\Phi},
\end{align}
\end{subequations}
where the combined integration measure $d\tilde{\Phi}$ absorbs the exact reference weights and the collision Jacobian:
\begin{equation}
    d\tilde{\Phi} = \maxwellianref(\vmag) \intweightref(\Ipre) \maxwellianref(\wmag) \intweightref(\Jpre) \phipart(\partint) \psipart(\partkin) (1-\partkin) \partkin^{1/2} \, d\partkin \, d\partint \, d\scatterdir \, d\Jpre \, d\wpre \, d\Ipre \, d\vpre.
\end{equation}

If the basis is truncated at maximum degrees $\kmax$, $\lmax$, and $\imax$, the total number of degrees of freedom scales as $\dof \propto (\kmax+1) (\imax+1) (\lmax+1)^2$. The resulting dense macroscopic tensor contains $\mathcal{O}(\kmax^3 \imax^3 \lmax^6)$ elements. This cubic growth motivates the exact dimensional reduction of Section~\ref{sec:reduction}.
\section{Dimensional Reduction of the Polyatomic Process}
\label{sec:reduction}

This section reduces the twelve-dimensional polyatomic collision integral to a nine-dimensional kinematic core. Integrating out the spatial orientation over the $\mathrm{SO}(3)$ rotation group decouples the intrinsic scattering and internal energy exchange from the macroscopic angular geometry, extending the Wigner--Eckart factorization of \citet{hiemstra2026wigner} to polyatomic gases.

\subsection{Geometric Symmetries and Dimensionality Reduction}

Evaluating the polyatomic collision operator requires integrating over a twelve-dimensional space defined by the pre-collision velocities $\vpre, \wpre \in \mathbb{R}^3$, the scattering direction $\scatterdir \in \mathbb{S}^2$, the microscopic internal energies $\Ipre, \Jpre \in [0, \infty)$, and the partition parameters $\partkin, \partint \in [0,1]$. The binary scattering kinematics, however, possess no preferred spatial orientation: the transition probability rate $B$ depends only on the relative speed $\urel$, the polar deflection angle $\defangle$, and the internal energy scalars.

\begin{figure}[t]
\centering
\includegraphics[width=\textwidth]{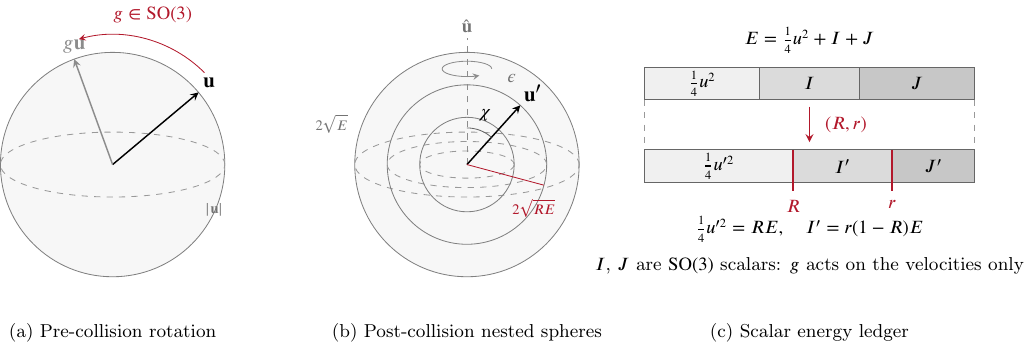}
\caption{$\mathrm{SO}(3)$ symmetry of the polyatomic collision. (a)~A rotation $g \in \mathrm{SO}(3)$ acts on the velocity frame as a whole, taking the pair $(\vpre, \wpre)$ to $(g\vpre, g\wpre)$; the relative velocity $\mathbf{u} = \vpre - \wpre$ inherits the rotation, and the collision map is equivariant under this action. A valid pre-collision configuration therefore maps to a valid configuration, with $\mathbf{u}$ and its rotated copy $g\mathbf{u}$ on the same sphere of radius $|\mathbf{u}|$. (b)~Post-collision, the kinetic partition $\partkin$ selects the scattering sphere $|\mathbf{u}'| = \sqrt{2\partkin\Etot}$ from the continuum $0 \leq \urel' \leq \sqrt{2\Etot}$ (a nested family replacing the single monatomic sphere $|\mathbf{u}'| = |\mathbf{u}|$), and the angles $(\defangle, \aziangle)$ select the point on it. The orbit structure is preserved sphere by sphere. (c)~The internal energies are $\mathrm{SO}(3)$ scalars: rotations act on the velocities only, so the two-stage exchange $\Etot = \tfrac{1}{2}\urel^2 + \Ipre + \Jpre \to \tfrac{1}{2}(\urel')^2 + \Ipost + \Jpost$ commutes with every $g$.}
\label{fig:so3_symmetry}
\end{figure}

Exploiting this rotational invariance eliminates three spatial degrees of freedom. The scalar internal variables ($\Ipre$, $\Jpre$, $\partkin$, $\partint$) are invariant under any rigid-body rotation, and the collision map is equivariant under the action of $g \in \mathrm{SO}(3)$ on the velocity pair: rotating the configuration rotates the outcome without changing the transition rate (Figure~\ref{fig:so3_symmetry}). We establish a privileged coordinate system by rigidly rotating the absolute laboratory frame to align with the colliding pair:
\begin{enumerate}
    \item \textbf{Alignment of the target particle:} the local $z$-axis is aligned with the target velocity, $\vpre = \vmag \hat{\mathbf{z}}$, factoring out the polar and azimuthal orientation of the collision event, eliminating two spatial degrees of freedom.
    \item \textbf{Alignment of the incident particle:} by cylindrical symmetry around the $z$-axis, the frame is rotated about $\vpre$ until the incident velocity $\wpre$ lies in the $x$-$z$ plane, integrating out its azimuthal degree of freedom.
\end{enumerate}

Having factored out the three-dimensional rigid-body rotation of the collision pair, the intrinsic geometry of the polyatomic collision is described by nine independent scalar variables:
\begin{itemize}
    \item $\vmag, \wmag$: The speeds of the target and incident particles.
    \item $\incangle$: The polar incidence angle between $\vpre$ and $\wpre$ (constrained to the $x$-$z$ plane).
    \item $\defangle, \aziangle$: The polar deflection and azimuthal scattering angles of $\scatterdir$ relative to $\vrelpre$.
    \item $\Ipre, \Jpre$: The pre-collision internal energies of the target and incident particles.
    \item $\partkin, \partint$: The kinetic and internal energy partition parameters.
\end{itemize}

In the monatomic setting, the same integration over the rigid-body rotation decouples the macroscopic angular geometry from a five-dimensional spatial kinematic core \citep{hiemstra2026wigner}. The polyatomic reduction differs in two respects: the invariant core expands to nine dimensions to carry the internal energy variables ($\Ipre, \Jpre$) and the partition parameters ($\partkin, \partint$), and the spectral basis is augmented with the internal polynomials $\intbasis_i$. The polyatomic complexity therefore resides in the rotationally invariant scattering dynamics, and the macroscopic angular decoupling is identical to the monatomic case.

\subsection{Exact Factorization via Integration over the Rotation Group}
Any collision configuration in the laboratory frame can be generated by applying a three-dimensional rigid rotation $g \in \mathrm{SO}(3)$ to the privileged reference configuration. Expressing the spatial integration measure in spherical coordinates isolates the invariant Haar measure $d g = \sin \theta \, d\theta \, d\phi \, d\psi$ over $\mathrm{SO}(3)$, yielding the geometric relation
\begin{equation}
    d\vpre \, d\wpre \, d\scatterdir = d g \, d\Phi_\mathrm{5D},
\end{equation}
where $d\Phi_\mathrm{5D} = \vmag^2 \wmag^2 \sin\incangle \sin\defangle \, d\vmag \, d\wmag \, d\incangle \, d\defangle \, d\aziangle$ encompasses the five-dimensional spatial kinematic core. The total differential measure from the Cartesian formulation thus strictly decouples into the macroscopic angular measure and a nine-dimensional intrinsic measure $d\Phi_\mathrm{9D} = d\Phi_\mathrm{5D} \, d\Psi_\mathrm{4D}$, with the internal measure defined as $d\Psi_\mathrm{4D} = d\partkin \, d\partint \, d\Jpre \, d\Ipre$.

When the continuous collision geometry is rotated by $g$, the scalar quantities remain invariant: the relative speed $\urel$, the transition kernel $B$, the radial Laguerre polynomials $\radbasis_{k,l}$, and the internal energy polynomials $\intbasis_i$. Only the spherical harmonics transform linearly according to the Wigner D-matrices $\wignerD_{mm'}^l(g)$ \citep{Varshalovich:1988}:
\begin{equation}
    \angbasis_{l,m}(\vpreunit) = \sum_{m'=-l}^l \wignerD_{mm'}^l(g) \angbasis_{l,m'}(\vpreunit_0),
\end{equation}
where the subscript $0$ denotes the directional vector evaluated in the privileged frame.

Substituting the transformed spherical harmonics decouples the integration over the rigid-body rotation ($d g$) from the invariant nine-dimensional core ($d\Phi_\mathrm{9D}$). By the representation identity for the product of three Wigner D-matrices over $\mathrm{SO}(3)$ \citep{Varshalovich:1988}, the rotational integral evaluates analytically to a product of Wigner 3-j symbols. The collision tensor then separates into a sparse macroscopic angular tensor $\Gtens$ and a dense, rotationally invariant physical tensor $\Rtens^+$:
\begin{equation}
    \Ctens_{\state_1, \state_2, \state_3}^+ = \Gtens_{m_1 m_2 m_3}^{l_1 l_2 l_3} \Rtens_{k_1 k_2 k_3, i_1 i_2 i_3}^{+ l_1 l_2 l_3}.
\end{equation}

The Gaunt tensor $\Gtens$ encapsulates the macroscopic angular conservation laws and remains identical to the monatomic formulation:
\begin{equation}
    \Gtens_{m_1 m_2 m_3}^{l_1 l_2 l_3} = 8\pi^2 \begin{pmatrix} l_1 & l_2 & l_3 \\ m_1 & m_2 & m_3 \end{pmatrix}.
\end{equation}
The $3$-$j$ symbols are evaluated in closed form, and the selection rules determine exactly which transitions survive; the generator stores the surviving entries in coordinate format.

The reduced physical tensor $\Rtens^+$ absorbs the internal energy dependencies, the collision cross-section, and the intrinsic scattering dynamics, evaluated over the invariant 9D core:
\begin{equation}
\begin{aligned}
    \Rtens_{k_1 k_2 k_3, i_1 i_2 i_3}^{+ l_1 l_2 l_3} &= \int_{\mathrm{9D}} B(\urel, \cos\defangle, \Ipre, \Jpre, \partkin, \partint) P_{\mathrm{gain}}^{l_1 l_2 l_3}(\vpostunit_0, \vpreunit_0, \wpreunit_0) \\
    &\quad \times \radbasis_{k_1, l_1}(\vmag') \intbasis_{i_1}(\Ipost) \radbasis_{k_2, l_2}(\vmag) \intbasis_{i_2}(\Ipre) \radbasis_{k_3, l_3}(\wmag) \intbasis_{i_3}(\Jpre) \, d\tilde{\Phi}_\mathrm{9D}.
\end{aligned}
\end{equation}

Here, $P_{\mathrm{gain}}^{l_1 l_2 l_3}(\vpostunit_0, \vpreunit_0, \wpreunit_0)$ is the invariant geometric filter, obtained by evaluating the three local spherical harmonics in the privileged frame and summing over the invariant local projection indices $m'$. Because the internal variables leave the collision vectors $\vpostunit$, $\vpreunit$, $\wpreunit$ untouched, the filter is identical to that of a monatomic gas. Its explicit derivation and Clebsch-Gordan expansion are given in \citet{hiemstra2026wigner}. The combined integration measure $d\tilde{\Phi}_\mathrm{9D}$ incorporates the nine-dimensional Jacobian alongside the associated Maxwellian and internal energy reference weights:
\begin{equation}
    d\tilde{\Phi}_\mathrm{9D} = \maxwellianref(\vmag) \intweightref(\Ipre) \maxwellianref(\wmag) \intweightref(\Jpre) \phipart(\partint) \psipart(\partkin) (1-\partkin) \partkin^{1/2} \, d\Phi_\mathrm{9D}.
\end{equation}

\subsection{Incorporation of the Loss Term}

The loss component $\Ctens^-$ undergoes the same rotational decoupling. Because the loss operator evaluates only pre-collision states, the associated spherical harmonics align with the pole of the privileged frame, collapsing the geometric loss filter $P_{\mathrm{loss}}^{l_1 l_2 l_3}$ into a zonal mode independent of the scattering angles \citep{hiemstra2026wigner}.

Since gain and loss share the Gaunt tensor $\Gtens$ and the integration measure $d\tilde{\Phi}_\mathrm{9D}$, they are merged at the continuous level into the net physical tensor $\Rtens = \Rtens^+ - \Rtens^-$:
\begin{equation}
\begin{aligned}
    \Rtens_{k_1 k_2 k_3, i_1 i_2 i_3}^{l_1 l_2 l_3} &= \int_{\mathrm{9D}} B(\urel, \cos\defangle, \Ipre, \Jpre, \partkin, \partint) \radbasis_{k_2, l_2}(\vmag) \intbasis_{i_2}(\Ipre) \radbasis_{k_3, l_3}(\wmag) \intbasis_{i_3}(\Jpre) \\
    &\quad \times \left[ P_{\mathrm{gain}}^{l_1 l_2 l_3}(\vpostunit_0, \vpreunit_0, \wpreunit_0) \radbasis_{k_1, l_1}(\vmag') \intbasis_{i_1}(\Ipost) \right. \\
    &\quad \left. - P_{\mathrm{loss}}^{l_1 l_2 l_3}(\vpreunit_0, \vpreunit_0, \wpreunit_0) \radbasis_{k_1, l_1}(\vmag) \intbasis_{i_1}(\Ipre) \right] \, d\tilde{\Phi}_\mathrm{9D}.
\end{aligned}
\end{equation}

This exact factorization reduces the $\mathcal{O}(\kmax^3 \imax^3 \lmax^6)$ complexity of the Cartesian collision tensor into a sparse $\mathcal{O}(\lmax^5)$ macroscopic assembly coupled with the evaluation of the dense $\mathcal{O}(\kmax^3 \imax^3 \lmax^3)$ net physical tensor $\Rtens$. The $\mathcal{O}(\lmax^5)$ count is the number of Gaunt transitions admitted by the triangle, parity, and azimuthal selection rules: the $\mathcal{O}(\lmax^3)$ polar triplets carry an $\mathcal{O}(\lmax^2)$ azimuthal multiplier, one order below $\lmax^3$ because the Kronecker delta on $m_3$ removes a summation degree of freedom. The count is unchanged by the internal states, and is derived in Section~4.6 of \citet{hiemstra2026wigner}.
\section{Numerical Implementation and Algorithm Design}
\label{sec:numerical}

This section constructs the numerical evaluation of the nine-dimensional integral defining the net physical tensor $\Rtens$. Tensorized Gauss-type quadrature over the invariant kinematic core reduces the continuous transition rates and internal energy exchanges to a discrete, precomputable array. The implementation is available as an open-source MATLAB library with C++ extensions \citep{kessler_zenodo_2026}.

The concrete transition kernel targeted throughout this section is the extended direct simulation Monte Carlo (DSMC) model of \citet{djordjic:2023}. It is a convex split of a \emph{non-frozen} channel (full Borgnakke--Larsen internal energy redistribution) and a \emph{frozen} channel (elastic, internal energy preserving), each modulated by fractional powers of the internal energy fractions:
\begin{equation}
\begin{aligned}
    B = \urel^{\vhs} \, \partkin^{\vhs/2} \Big[\, & \inelprob \, \Kmono \, \Big( 1 + \modamp \Big[ \Big( \partint (1-\partkin) \tfrac{\Ipre}{\Etot} \Big)^{\modexp} + \Big( (1-\partint)(1-\partkin) \tfrac{\Jpre}{\Etot} \Big)^{\modexp} \Big] \Big) \\
    & + (1-\inelprob) \, \delta(\partkin - \bar{\partkin}) \, \delta(\partint - \bar{\partint}) \Big( 1 + \modampf \Big[ \Big( \tfrac{\Ipre}{\Etot} \Big)^{\modexpf} + \Big( \tfrac{\Jpre}{\Etot} \Big)^{\modexpf} \Big] \Big) \Big],
\end{aligned}
\label{eq:dsmc_kernel}
\end{equation}
where $\vhs = 2(1 - s_{\mathrm{visc}})$ is the variable-hard-sphere (VHS) exponent \citep{bird:1994}, with $s_{\mathrm{visc}}$ the temperature exponent of the viscosity law, and $\inelprob \in [0,1]$ is the inelastic collision probability. The remaining factors of \eqref{eq:dsmc_kernel} are:
\begin{itemize}
    \item \textbf{VHS collision rate ($\urel^{\vhs}$):} multiplies both channels (energy weight, Section~\ref{sec:quadfamilies}).
    \item \textbf{Post-collision partition density ($\partkin^{\vhs/2}$):} the density of the post-collision translational energy fraction, multiplying both channels. On the non-frozen channel it is the Borgnakke--Larsen partition density (kinetic-partition weight, Section~\ref{sec:quadfamilies}). On the frozen channel the Dirac factors below fix $\partkin$ at its elastic value, so the same factor is not integrated but evaluated pointwise as the weight $\bar{\partkin}^{\vhs/2} = (\urel^2/2\Etot)^{\vhs/2}$ (auxiliary integral, Section~\ref{sec:laplace}).
    \item \textbf{Measure compensation ($\Kmono = 2\Gamma(\dofpoly + \tfrac{3}{2}) / [\sqrt{\pi} \, \Gamma(\dofpoly/2)^2]$):} constant, non-frozen channel only, with $\dofpoly = 2\polyalpha + 2$ the number of internal degrees of freedom (Section~\ref{sec:background}). It compensates the Borgnakke--Larsen partition measure so that the channel's elastic-limit rate equals the plain VHS rate.
    \item \textbf{Elastic pinning ($\delta(\partkin - \bar{\partkin}) \, \delta(\partint - \bar{\partint})$):} pins the frozen channel at $\bar{\partkin} = \urel^2/(2\Etot)$, $\bar{\partint} = \Ipre/(\Ipre + \Jpre)$, reproducing the elastic collision $\urel' = \urel$, $\Ipost = \Ipre$, $\Jpost = \Jpre$ (partition collapse, Appendix~\ref{app:assembly}).
    \item \textbf{Non-frozen modulation ($\modamp$, $\modexp$):} internal-energy modulation with amplitude $\modamp \geq -\tfrac{1}{2}$ (kernel positivity) and exponent $\modexp \geq 0$ fitted to experimental transport data (auxiliary integral, Section~\ref{sec:laplace}).
    \item \textbf{Frozen modulation ($\modampf$, $\modexpf$):} the frozen-channel counterpart, $\modampf \geq -\tfrac{1}{2}$, $\modexpf \geq 0$ (auxiliary integral, Section~\ref{sec:laplace}). Our exponent convention $(\modexp, \modexpf)$ corresponds to $(\hat{\zeta}/2, \hat{\zeta}_f)$ in the notation of \citet{djordjic:2023}.
\end{itemize}
The unweighted equilibrium averages of the two channels' partition factors are identical, because the pinned value $\bar{\partkin}$ and the Borgnakke--Larsen partition $\partkin$ both follow the $\mathrm{Beta}(3/2, \dofpoly)$ law. The two channel rates nevertheless differ, because the pinned factor is evaluated pointwise and is positively correlated with the kernel's own $\urel^{\vhs}$. For the base kernels of the three production gases of Section~\ref{sec:res_transport} the frozen channel's equilibrium rate is $1.033$ (N$_2$), $1.033$ (CO), and $1.041$ (H$_2$) times the non-frozen rate, and the modulation amplitudes reshape both rates further. This differs from the textbook DSMC picture, in which $\inelprob$ selects the inelastic treatment of a collision at a fixed total collision rate. Here $\inelprob$ sets the inelastic fraction of collisions, and the total rate depends on it.

Equation~\eqref{eq:dsmc_kernel} is Eq.~(43) of \citet{djordjic:2023}, the extended model from which their transport relations and fitted parameters are derived. The form quoted in their DSMC comparison, their Eq.~(53), carries no such factor on the frozen channel and is therefore not equivalent to it for $\vhs > 0$. We adopt one normalization convention with \eqref{eq:dsmc_kernel}: the Dirac factors of the frozen channel are normalized against the weighted partition measure, Eq.~\eqref{eq:frozen_normalization} of Appendix~\ref{app:assembly}, which removes the factor $\phipart(\bar{\partint}) \psipart(\bar{\partkin}) (1 - \bar{\partkin}) \bar{\partkin}^{1/2}$ that a literal reading of the deltas would leave. The frozen-limit agreement with their closed-form Eq.~(42), reported in Section~\ref{sec:res_transport_pred}, is the evidence that their evaluation used the same convention.

The discrete infrastructure of the monatomic implementation carries over unchanged, with the composite index $(k, l, m, i)$ replacing the monatomic $(k, l, m)$: the real-valued index mapping with its triangle, parity, and azimuthal selection rules, the coordinate (COO) memory layout, and the storage and algorithmic complexity analysis (Sections~4.1--4.2 and 4.6--4.7 of \citet{hiemstra2026wigner}). The sum-factorized offline assembly specific to the polyatomic build is organized in Appendix~\ref{app:assembly}. The remainder concentrates on the polyatomic ingredients: the quadrature architecture for the nine-dimensional core (Section~\ref{sec:quadarch}), the exact conservation embedding (Section~\ref{sec:conservation}), and the contraction over the enlarged radial--internal mode space (Section~\ref{sec:contraction}).

\subsection{Singular Quadrature Architecture}
\label{sec:quadarch}

Mirroring the monatomic architecture \citep{hiemstra2026wigner}, the physical tensor is evaluated on a tensor-product grid of mapped one-dimensional rules, one for each of the ten integration directions.

\subsubsection{Quadrature Formulation and Coordinate Transformations}
\label{sec:quadfamilies}

Our design principle, inherited from the monatomic formulation \citep{hiemstra2026wigner}, is composition. Every integration direction is assigned a transformation and a quadrature rule matched to its integrand class: polynomial exactness against the native weight, singularity-removing coordinate transformations, pointwise entire factors, or a geometrically convergent auxiliary rule. For the full kernel family of \eqref{eq:dsmc_kernel} every direction is then spectral, with a single scoped exception: the fractional endpoint power $\urel^{\vhs}$ left at the Duffy vertex, whose algebraic tail enters four to five orders below working accuracy (Section~\ref{sec:res_convergence}). The internal node counts are consequently exactness bounds rather than asymptotic estimates: the production path applies no padding to the internal axes, a choice verified immaterial to $2 \times 10^{-14}$.

Throughout, we work in the reference units of Section~\ref{sec:basis}, in which the Maxwellian weight is $e^{-\vmag^2}$ and the internal energies carry the canonical weight $\Ipre^{\polyalpha} e^{-\Ipre}$. Each of the ten integration directions is then treated as follows:
\begin{itemize}
    \item \textbf{Center-of-mass energy ($\cmenergy \in [0, \infty)$):} A $45^\circ$ rotation of the speed plane into center-of-mass and relative coordinates, composed with the square map $\cmenergy = \cmvar^2$, aligns the energy axis with the joint Maxwellian decay (Section~\ref{sec:kinematics}). The weight $W(\cmenergy) = \cmenergy^{\vhs/2} e^{-\cmenergy}$ collects the reference Maxwellians and the energy power of the kernel factor $\urel^{\vhs}$. \textbf{Generalized Gauss--Laguerre} quadrature with parameter $\alpha = \vhs/2$ is exact for the polynomial radial states.
    \item \textbf{Duffy patch pair ($\fracrel$ or $\halfangle$ outer, $\duffy \in [0,1]$ inner):} The fractional relative coordinate $\fracrel = \srel/\cmvar$ and the incidence half-angle $\halfangle = \sin(\incangle/2)$ parameterize the relative-motion square, on which the relative speed carries a cone singularity at the origin. A 2D Duffy transformation \citep{duffy:1982} splits the square into two triangular patches ($\halfangle = \fracrel \duffy$ and $\fracrel = \halfangle \duffy$), factoring the singularity into an explicit radial arm absorbed by the measure (Section~\ref{sec:kinematics}). On each patch the integrand is smooth, and \textbf{Gauss--Legendre} rules converge spectrally \citep{hiemstra2026wigner}.
    \item \textbf{Polar deflection ($\cos\defangle \in [-1,1]$):} The geometric filters and the reconstructed post-collision states are polynomial in $\cos\defangle$, so \textbf{Gauss--Legendre} quadrature is exact. The symmetry of the rule also participates in canceling odd powers of the post-collision speed (Appendix~\ref{app:assembly}).
    \item \textbf{Azimuthal angle ($\aziangle \in [0, 2\pi)$):} The integrand is a trigonometric polynomial of bounded degree on the periodic circle, for which the uniform \textbf{Trapezoidal} rule is spectrally exact.
    \item \textbf{Internal energies ($\Ipre, \Jpre \in [0, \infty)$):} The integration incorporates the reference equilibrium distributions. The weight $W(\Ipre) = \Ipre^\polyalpha e^{-\Ipre}$ is exactly the canonical form for \textbf{Generalized Gauss--Laguerre} quadrature with parameter $\alpha = \polyalpha$, exact for the polynomial internal basis.
    \item \textbf{Internal partition parameter ($\partint \in [0, 1]$):} The weight is dictated by the partition function $\phipart(\partint) = \partint^\polyalpha (1-\partint)^\polyalpha$. Mapping the domain to $[-1,1]$ recovers the standard Jacobi weight. We evaluate via \textbf{Gauss--Jacobi} quadrature with parameters $\alpha = \polyalpha$ and $\beta = \polyalpha$.
    \item \textbf{Kinetic partition parameter ($\partkin \in [0, 1]$):} The weight combines the partition function $\psipart(\partkin)$, the collision Jacobian, and the kernel factor $\partkin^{\vhs/2}$ of \eqref{eq:dsmc_kernel}, yielding $W(\partkin) = \partkin^{(1+\vhs)/2} (1-\partkin)^{2\polyalpha+1}$, evaluated via \textbf{Gauss--Jacobi} quadrature with parameters $\alpha = 2\polyalpha+1$ and $\beta = (1+\vhs)/2$.
    \item \textbf{Auxiliary Laplace direction ($\auxmap \in [0, 1]$):} Required for the fractional modulation factors of \eqref{eq:dsmc_kernel} and for the frozen channel's pinned partition factor, which couple the kinetic and internal energies and break the product structure above. A Laplace representation of the coupled power restores the product structure at the cost of a tenth integration direction, carrying the Jacobi weight $\auxmap^{\modexp-1}(1-\auxmap)^{2\polyalpha+1}$ in the non-frozen case (the frozen-channel parameters are given in Section~\ref{sec:laplace}) and evaluated by \textbf{Gauss--Jacobi} quadrature with $\numaux$ nodes. Its convergence is observed to be geometric rather than governed by a polynomial exactness bound (Section~\ref{sec:laplace}).
\end{itemize}

In the monatomic truncation $\imax = 0$ the four internal directions degenerate to a single-node Dirac-collapse rule: the internal energies are pinned at $\Ipre = \Jpre = 0$ with the full measure mass $\Gamma(\polyalpha + 1)$ retained, and the kinetic partition at $\partkin = 1$. With this rule the polyatomic build recovers the monatomic operator exactly (verified in Section~\ref{sec:res_mono}).

\begin{figure}[t]
\centering
\includegraphics[width=\textwidth]{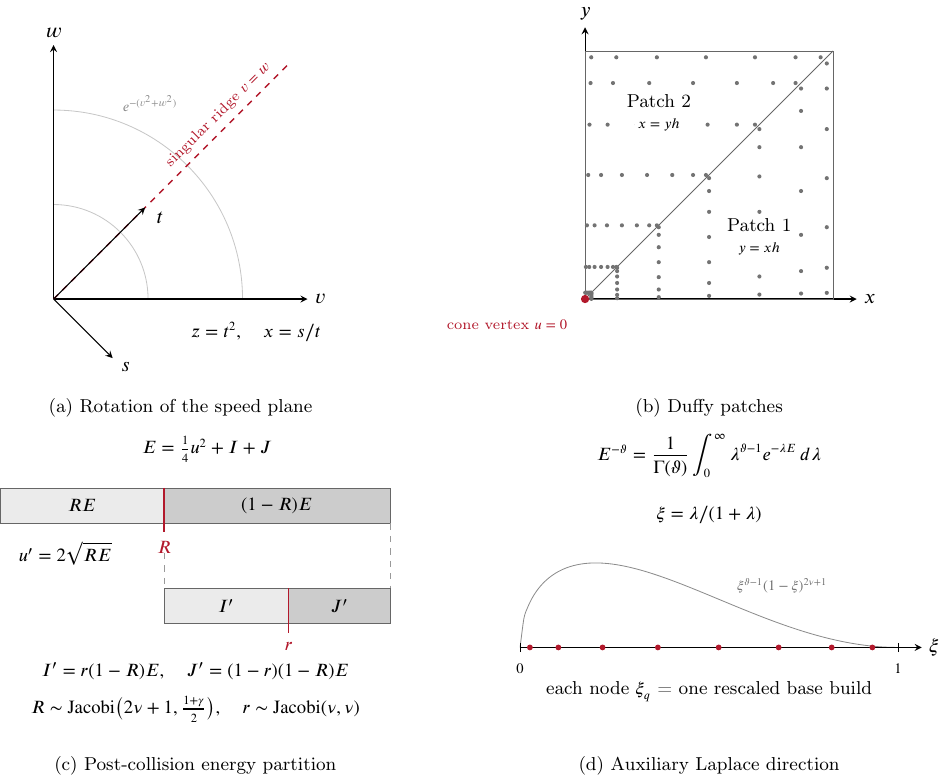}
\caption{Geometry of the quadrature construction. (a)~The $45^\circ$ rotation of the speed plane $(\vmag, \wmag)$ into the coordinates $\cmvar = (\vmag+\wmag)/\sqrt{2}$, $\srel = (\vmag-\wmag)/\sqrt{2}$ aligns the singular ridge $\vmag = \wmag$ of the kernel with a coordinate axis; the energy map $\cmenergy = \cmvar^2$ and the fractional coordinate $\fracrel = \srel/\cmvar$ follow. (b)~The Duffy transformation splits the relative-motion square $(\fracrel, \halfangle)$ into two triangular patches meeting at the cone vertex $\urel = 0$. The mapped tensor-product Gauss nodes (dots) accumulate toward the vertex, and the patch maps absorb the cone factor into the measure. (c)~The Borgnakke--Larsen partition splits the total collision energy $\Etot$ in two stages: the kinetic partition $\partkin$ fixes the post-collision relative speed $\urel' = \sqrt{2 \partkin \Etot}$, and the internal partition $\partint$ divides the remainder between $\Ipost$ and $\Jpost$, each stage integrated with the Gauss--Jacobi weight indicated. (d)~The auxiliary Laplace direction maps $\auxvar \in [0, \infty)$ to $\auxmap = \auxvar/(1+\auxvar) \in [0,1]$, carrying the Jacobi weight $\auxmap^{\modexp-1}(1-\auxmap)^{2\polyalpha+1}$; each Gauss--Jacobi node amounts to two rescaled builds of the base tensor, one for each active internal axis (Section~\ref{sec:laplace}).}
\label{fig:quadrature_geometry}
\end{figure}

Figure~\ref{fig:quadrature_geometry} summarizes the geometry of these transformations.

Because the asymmetric micro-reversibility quotient cancels analytically in the weak form (Section~\ref{sec:background}), the integrand of every direction above is smooth and bounded, and no truncation is required at the partition endpoints $\partint \in \{0,1\}$ or $\partkin = 1$.

The minimum node counts follow from the algebraic degrees dictated by the truncation limits $\kmax$, $\lmax$, and $\imax$. The five spatial coordinates inherit the monatomic bounds and the 10--20 point padding of the coupled coordinates $(\cmenergy, \fracrel, \halfangle, \duffy)$ \citep{hiemstra2026wigner}, while the four internal axes run at their exact-by-degree counts without padding (Table~\ref{tab:quadrature_bounds}).

\begin{table}[ht]
\centering
\caption{Summary of the mapped coordinates, weight functions, and minimum quadrature bounds required to evaluate the 9D polyatomic collision core. The spatial coordinates $(\cmenergy, \fracrel, \halfangle, \duffy, \cos\defangle, \aziangle)$ rely on the singularity-resolving architecture of the monatomic baseline \citep{hiemstra2026wigner}. The bounds are dictated by the radial ($\kmax$), angular ($\lmax$), and internal ($\imax$) spectral limits. $\vhs$ denotes the VHS exponent of the interaction potential, and $\polyalpha = \dofpoly/2 - 1$ denotes the internal-Laguerre parameter. The auxiliary direction $\auxmap$ applies to energy-modulated kernels and to the frozen channel's pinned partition factor (the non-frozen weight is shown; the frozen-channel parameters are given in Section~\ref{sec:laplace}). Its rule converges geometrically rather than through polynomial exactness.}
\label{tab:quadrature_bounds}
\resizebox{\textwidth}{!}{
\begin{tabular}{@{}llllc@{}}
\toprule
\textbf{Coordinate} & \textbf{Active Domain} & \textbf{Quadrature Rule} & \textbf{Weight Function} & \textbf{Minimum Points ($N$)} \\ \midrule
Energy ($\cmenergy$) & Both Patches & Gen. Gauss--Laguerre & $W(\cmenergy) = \cmenergy^{\vhs/2}e^{-\cmenergy}$ & $\lceil \frac{3\kmax + 1.5\lmax + 3}{2} \rceil$ \\
Outer Frac. ($\fracrel$) & Patch 1 ($\fracrel > \halfangle$) & Gauss--Legendre & $W(\fracrel) = 1$ & $4\kmax + 3\lmax + 4$ \\
Outer Angle ($\halfangle$) & Patch 2 ($\halfangle > \fracrel$) & Gauss--Legendre & $W(\halfangle) = 1$ & $4\kmax + 3\lmax + 4$ \\
Inner Duffy ($\duffy$) & Patch 1 ($\halfangle = \fracrel \duffy$) & Gauss--Legendre & $W(\duffy) = 1$ & $\kmax + \lceil 1.5\lmax \rceil + 1$ \\
Inner Duffy ($\duffy$) & Patch 2 ($\fracrel = \halfangle \duffy$) & Gauss--Legendre & $W(\duffy) = 1$ & $3\kmax + \lfloor 1.5\lmax \rfloor + 3$ \\
Deflection ($\cos\defangle$) & Both Patches & Gauss--Legendre & $W(\cos\defangle) = 1$ & $\kmax + \lceil 0.5\lmax \rceil + 1$ \\
Azimuthal ($\aziangle$) & Both Patches & Trapezoidal & $W(\aziangle) = 1$ & $2\kmax + \lmax + 1$ \\ \midrule
Target Internal ($\Ipre$) & $[0, \infty)$ & Gen. Gauss--Laguerre & $W(\Ipre) = \Ipre^\polyalpha e^{-\Ipre}$ & $\lceil \frac{\kmax + 2\imax + 1}{2} \rceil$ \\
Incident Internal ($\Jpre$) & $[0, \infty)$ & Gen. Gauss--Laguerre & $W(\Jpre) = \Jpre^\polyalpha e^{-\Jpre}$ & $\lceil \frac{\kmax + 2\imax + 1}{2} \rceil$ \\
Internal Partition ($\partint$) & $[0, 1]$ & Gauss--Jacobi & $W(\partint) = \partint^\polyalpha(1-\partint)^\polyalpha$ & $\imax + 1$ \\
Kinetic Partition ($\partkin$) & $[0, 1]$ & Gauss--Jacobi & $W(\partkin) = \partkin^{(1+\vhs)/2}(1-\partkin)^{2\polyalpha+1}$ & $\lceil \frac{\kmax + 2\imax + 1}{2} \rceil$ \\ \midrule
Auxiliary ($\auxmap$) & $[0, 1]$ & Gauss--Jacobi & $W(\auxmap) = \auxmap^{\modexp-1}(1-\auxmap)^{2\polyalpha+1}$ & geometric; $\numaux = 24$ for $10^{-12}$, $40$ for roundoff (Sec.~\ref{sec:laplace}) \\ \bottomrule
\end{tabular}
}
\end{table}

\subsubsection{Kinematic Mapping at the Quadrature Nodes}
\label{sec:kinematics}

The tensorized rules place nodes in the mapped coordinates, whereas the integrand of the net physical tensor is expressed in the physical collision variables. The two are connected by a chain of closed-form reconstructions evaluated at every node. The $45^\circ$ rotation into center-of-mass and relative speed coordinates, $\cmvar = (\vmag + \wmag)/\sqrt{2}$ and $\srel = (\vmag - \wmag)/\sqrt{2}$, composed with the energy map $\cmenergy = \cmvar^2$ and the fractional coordinate $\fracrel = \srel/\cmvar \in [0,1]$, inverts explicitly to the pre-collision speeds
\begin{equation}
    \vmag(\cmenergy, \fracrel) = \sqrt{\tfrac{\cmenergy}{2}} \, (1 + \fracrel), \qquad \wmag(\cmenergy, \fracrel) = \sqrt{\tfrac{\cmenergy}{2}} \, (1 - \fracrel).
    \label{eq:speed_reconstruction}
\end{equation}
The incidence angle enters through the half-angle parameter $\halfangle = \sin(\incangle/2)$, so that $\cos\incangle = 1 - 2\halfangle^2$ and the pre-collision relative speed evaluates to $\urel = \sqrt{(\vmag - \wmag)^2 + 4 \vmag \wmag \halfangle^2} = \sqrt{2\cmenergy} \sqrt{\fracrel^2 + (1 - \fracrel^2)\halfangle^2}$. Under the two Duffy patch maps that resolve the cone singularity of this square root at $\fracrel = \halfangle = 0$, the relative speed takes the per-patch closed forms
\begin{equation}
    \urel = \sqrt{2\cmenergy} \, \fracrel \sqrt{1 + (1 - \fracrel^2) \duffy^2} \quad (\halfangle = \fracrel \duffy), \qquad
    \urel = \sqrt{2\cmenergy} \, \halfangle \sqrt{\duffy^2 + 1 - \halfangle^2 \duffy^2} \quad (\fracrel = \halfangle \duffy),
    \label{eq:duffy_speeds}
\end{equation}
in which the radial arm of the singularity appears as the explicit prefactor $\fracrel$ (respectively $\halfangle$) absorbed by the patch measure in the assembly (Appendix~\ref{app:assembly}).

The Borgnakke--Larsen mechanics of Section~\ref{sec:background} then evaluate along the dependency chain
\begin{equation}
    \Etot = \tfrac{1}{2}\urel^2 + \Ipre + \Jpre = \cmenergy \big( \fracrel^2 + (1 - \fracrel^2)\halfangle^2 \big) + \Ipre + \Jpre, \qquad
    \urel' = \sqrt{2 \partkin \Etot}, \qquad
    \Ipost = \partint (1 - \partkin) \Etot,
    \label{eq:energy_chain}
\end{equation}
and the post-collision target state is reconstructed geometrically from $\vpost = \vcm + \tfrac{1}{2} \urel' \scatterdir$, where $\vcm = \tfrac{1}{2}(\vpre + \wpre)$ is the pair center-of-mass velocity and $\scatterdir(\defangle, \aziangle)$ is the scattering direction built in the privileged frame of Section~\ref{sec:reduction}:
\begin{equation}
    (\vmag')^2 = |\vcm|^2 + \tfrac{1}{4}(\urel')^2 + \urel' \, (\vcm \cdot \scatterdir), \qquad
    \vpostunit = \frac{\vcm + \tfrac{1}{2} \urel' \scatterdir}{\vmag'}.
    \label{eq:post_reconstruction}
\end{equation}
Both $|\vcm|$ and $\vcm \cdot \scatterdir$ are closed-form functions of the node coordinates in the privileged frame, identical to the monatomic reconstruction \citep{hiemstra2026wigner}. The chain exposes the functional dependencies that organize the sum-factorized assembly of Appendix~\ref{app:assembly}: $\Etot$ depends only on $(\cmenergy, \fracrel, \halfangle, \Ipre, \Jpre)$, the post-collision speed $\vmag'$ and direction $\vpostunit$ add the kinetic partition $\partkin$ and the scattering angles, and only $\Ipost$ depends on the internal partition $\partint$.

\subsubsection{Energy-Coupled Kernels and the Auxiliary Laplace Integral}
\label{sec:laplace}

The modulation factors of both channels augment the transition rate, up to the bounded powers of the partition variables accompanying them in \eqref{eq:dsmc_kernel}, by terms of the form
\begin{equation}
    B \mapsto B \left[ 1 + \modamp \left( \left( \frac{\Ipre}{\Etot} \right)^{\modexp} + \left( \frac{\Jpre}{\Etot} \right)^{\modexp} \right) \right],
    \label{eq:modulated_kernel}
\end{equation}
where the exponent $\modexp > 0$ is in general non-integer. We develop the treatment for the generic form \eqref{eq:modulated_kernel}; the frozen channel follows by the substitution $(\modamp, \modexp) \to (\modampf, \modexpf)$.

Written as $\Ipre^{\modexp} \Etot^{-\modexp}$, the modulation couples independent integration axes through the total collision energy $\Etot = \frac{1}{2}\urel^2 + \Ipre + \Jpre$. The branch point of $\Etot^{-\modexp}$ at $\Etot \to 0$ is not removable by any reparameterization of the product domain, and pointwise evaluation at the tensorized nodes reduces the rules to algebraic convergence.

We resolve the coupling by integration rather than reparameterization, through the Gamma-function (Laplace) identity
\begin{equation}
    \Etot^{-\modexp} = \frac{1}{\Gamma(\modexp)} \int_0^\infty \auxvar^{\modexp - 1} e^{-\auxvar \Etot} \, d\auxvar, \qquad \modexp > 0.
    \label{eq:laplace_identity}
\end{equation}
Because the total energy is additive with no cross terms, the auxiliary exponential factorizes across the independent axes,
\begin{equation}
    e^{-\auxvar \Etot} = e^{-\auxvar \urel^2 / 2} \, e^{-\auxvar \Ipre} \, e^{-\auxvar \Jpre},
    \label{eq:laplace_factorization}
\end{equation}
so that for every fixed $\auxvar$ the nine-dimensional integrand recovers the product structure required by the tensorized rules.

Each axis absorbs its factor without introducing new singularities. In the reference units, the fractional power $\Ipre^{\modexp}$ on the active internal axis shifts the generalized Laguerre parameter from $\polyalpha$ to $\polyalpha + \modexp$, while the exponential $e^{-\auxvar \Ipre}$ rescales the quadrature nodes by $\auxscale = (1+\auxvar)^{-1}$:
\begin{equation}
    \int_0^\infty \Ipre^{\polyalpha + \modexp} e^{-(1+\auxvar)\Ipre} \, g(\Ipre) \, d\Ipre = \auxscale^{\polyalpha + \modexp + 1} \int_0^\infty \tau^{\polyalpha + \modexp} e^{-\tau} g(\auxscale \tau) \, d\tau.
    \label{eq:node_rescaling}
\end{equation}
The spectator axis $\Jpre$ is treated identically with the unshifted parameter $\polyalpha$, and both axes remain spectral for every $\auxvar$. The velocity factor $e^{-\auxvar \urel^2/2}$ is entire and applied pointwise on the unchanged spatial architecture. The bounded partition prefactors accompanying \eqref{eq:modulated_kernel} fold into parameter-shifted Gauss--Jacobi weights. The factor $(1-\partkin)^{\modexp}$ shifts the kinetic-partition parameters to $(\alpha, \beta) = (2\polyalpha + 1 + \modexp, \, (1+\vhs)/2)$, and the factor $\partint^{\modexp}$ (respectively $(1-\partint)^{\modexp}$) shifts the internal-partition parameters to $(\polyalpha, \, \polyalpha + \modexp)$ for the $\Ipre$-active channel and $(\polyalpha + \modexp, \, \polyalpha)$ for the $\Jpre$-active channel. The remaining partition integrands stay polynomial and are therefore integrated exactly. The two modulation terms define an $\Ipre$-active and a $\Jpre$-active channel, which are summed.

The auxiliary integral itself is mapped to the unit interval via the substitution $\auxmap = \auxvar/(1+\auxvar)$, so that $\auxscale = 1 - \auxmap$. Collecting the node-scaling powers $\auxscale^{2\polyalpha + \modexp + 2}$ contributed by the two internal axes together with the transformed measure $\auxvar^{\modexp-1} \, d\auxvar$ leaves the canonical Jacobi weight on $[0,1]$:
\begin{equation}
    \Rtens_{\mathrm{corr}} = \frac{1}{\Gamma(\modexp)} \int_0^1 \auxmap^{\modexp - 1} (1-\auxmap)^{2\polyalpha + 1} \, \widehat{\Rtens}(\auxmap) \, d\auxmap,
    \label{eq:aux_jacobi}
\end{equation}
evaluated by Gauss--Jacobi quadrature with parameters $\alpha = 2\polyalpha + 1$ and $\beta = \modexp - 1$. Here $\widehat{\Rtens}(\auxmap)$ denotes one evaluation of the base physical tensor, the tensorized quadrature of Section~\ref{sec:quadfamilies} applied to the unmodulated kernel, with the internal nodes rescaled according to \eqref{eq:node_rescaling} and the velocity factor of \eqref{eq:laplace_factorization} applied pointwise. The auxiliary rule converges at approximately one decade per three to four nodes over $\numaux = 10$ to $24$, the range over which the rate is quoted, and flattens beyond it. At $\numaux = 16$ the error still runs from $3.6 \times 10^{-11}$ to $6.8 \times 10^{-10}$ across the branches. Double precision requires $\numaux \geq 40$.

Because the operator is linear in the modulation amplitudes, the extended physical tensor assembles as $\Rtens = \Rtens_{\mathrm{base}} + \modamp \, \Rtens_{\mathrm{corr}}$ (and analogously with $\modampf$ in the frozen channel), where $\Rtens_{\mathrm{base}}$ belongs to the base variable-hard-sphere Borgnakke--Larsen kernel recovered at $\modamp = \modampf = 0$. The validated base build is thus reused unchanged, and the correction amounts to $2\numaux$ further evaluations of the base tensor at the rescaled nodes with the shifted weights, one $\Ipre$-active and one $\Jpre$-active per auxiliary node. The auxiliary factor multiplies the gain and loss terms identically, so the exact embedding constructed next (Section~\ref{sec:conservation}) applies verbatim to the extended tensor.

The frozen channel's pinned partition factor $\bar{\partkin}^{\vhs/2}$ is treated by the same device, with the exponents rearranged by the split $\bar{\partkin}^{\vhs/2} = (\urel^2/2)^{\vhs/2} \, \Etot^{-\vhs/2}$. The translational factor combines with the kernel factor $\urel^{\vhs}$, so the frozen integrand carries $\urel^{2\vhs}$ and the energy axis runs the generalized Gauss--Laguerre rule with parameter $\vhs$ in place of $\vhs/2$. The coupled factor $\Etot^{-\vhs/2}$ passes through the identity \eqref{eq:laplace_identity} with exponent $p = \vhs/2$. On the frozen modulated branch the additional power $\Etot^{-\modexpf}$ combines with it into the single exponent $p = \vhs/2 + \modexpf$ (one auxiliary integral, not two), while the numerator $\Ipre^{\modexpf}$ shifts the active internal Laguerre parameter to $\polyalpha + \modexpf$ exactly as in \eqref{eq:node_rescaling}. Collecting the node-scaling powers as in \eqref{eq:aux_jacobi} leaves, on both branches, the Jacobi weight $\auxmap^{p-1}(1-\auxmap)^{2\polyalpha+1-\vhs/2}$, evaluated with parameters $\alpha = 2\polyalpha + 1 - \vhs/2$ and $\beta = p - 1$. The measured convergence matches the non-frozen rule: approximately one decade per three to four auxiliary nodes over $\numaux = 10$ to $24$. Measured in the relative $\ell_\infty$ norm of the physical tensor against an over-resolved $\numaux = 44$ reference, the error at $\numaux = 24$ is at most $8 \times 10^{-13}$ across the frozen base and modulated branches of all three gases, including the largest production exponent $p = 1.27$ (H$_2$: $\vhs = 0.607$, $\modexpf = 0.965$), and at most $3 \times 10^{-12}$ for the isolated modulated correction. Every branch reaches the roundoff floor at $\numaux = 40$.

The organization of these evaluations into an algorithm is detailed in Appendix~\ref{app:assembly}: a strict offline/online split in which the tensorized quadrature is sum-factorized along the dependency chain of Section~\ref{sec:kinematics} into nested one-dimensional stages, with the frozen channel as a separate elastic pass. The offline build costs $1 + 2\numaux$ base builds per kernel calibration and stores $\mathcal{O}(\kmax^3 \imax^3 \lmax^3)$ physical-tensor entries.

\subsection{Enforcement of Macroscopic Conservation}
\label{sec:conservation}

The five collision invariants (mass, momentum, and the total translational plus internal energy) must annihilate the discrete operator exactly. Because the macroscopic collision tensor of Section~\ref{sec:reduction} is linear in the physical tensor through the shared test index $(k_1, l_1, i_1)$, constraining rows of $\Rtens$ is equivalent to constraining rows of the assembled operator, at a fraction of the cost. Mass and momentum occupy the single test row $(k_1, i_1) = (0,0)$: this row is zeroed outright in the $l_1 = 0$ (mass) and $l_1 = 1$ (momentum) channels, exactly as in the monatomic embedding \citep{hiemstra2026wigner}.

The total energy spans two test rows. In the reference units of Section~\ref{sec:basis}, the Laguerre identities $L_1^{(1/2)}(\vmag^2) = \tfrac{3}{2} - \vmag^2$ and $L_1^{(\polyalpha)}(\Ipre) = (\polyalpha + 1) - \Ipre$ give the exact two-mode expansion of the energy invariant in the orthonormal reference basis,
\begin{equation}
    \vmag^2 + \Ipre = \Big( \tfrac{3}{2} + \polyalpha + 1 \Big) \basis_{(0,0,0,0)} - \sqrt{\tfrac{3}{2}} \, \basis_{(1,0,0,0)} - \sqrt{\polyalpha + 1} \, \basis_{(0,0,0,1)},
    \label{eq:energy_invariant}
\end{equation}
with the composite index $\state = (k, l, m, i)$. Energy conservation is therefore the single linear constraint $\sqrt{3/2} \cdot \Rtens_{(1,0)} + \sqrt{\polyalpha+1} \cdot \Rtens_{(0,1)} = 0$ on the row pair $\Rtens_{(k_1, i_1)}$ with $(k_1, i_1) \in \{(1,0), (0,1)\}$ in the $l_1 = 0$ channels, imposed by the orthogonal projection
\begin{equation}
    \begin{pmatrix} \Rtens_{(1,0)} \\ \Rtens_{(0,1)} \end{pmatrix}
    \;\leftarrow\;
    \big( \mathrm{Id} - \hat{e} \hat{e}^{\mathsf{T}} \big)
    \begin{pmatrix} \Rtens_{(1,0)} \\ \Rtens_{(0,1)} \end{pmatrix},
    \qquad
    \hat{e} \propto \begin{pmatrix} \sqrt{3/2} \\ \sqrt{\polyalpha + 1} \end{pmatrix}.
    \label{eq:energy_projection}
\end{equation}

The projection removes only the invariant combination. The orthogonal complement of the row pair carries the physical translational--internal energy exchange (Landau--Teller) mode and is left exactly as computed, so the relaxation physics of Section~\ref{sec:results} is unaffected by the embedding. In the monatomic limit $\imax = 0$ the internal row is absent, the constraint reduces to zeroing the $(k_1, i_1) = (1, 0)$ row alone, and the monatomic per-row embedding is recovered. The enforced linearized operator possesses exactly the five-dimensional null space spanned by the invariants. In the frozen-only limit $\inelprob = 0$ the internal energy is separately conserved, so every internal mode $i = 1, \ldots, \imax$ contributes a further invariant and the null space is $(5 + \imax)$-dimensional, which is seven-dimensional at the production truncation $\imax = 2$.

\subsection{Tensor Contraction Strategies}
\label{sec:contraction}

With the physical tensor integrated and conservation-corrected, the online operator evaluates the discrete factored summation
\begin{equation}
    Q_{(k_1 i_1), q_1} = \sum_{\{z \,\mid\, q_1^z = q_1\}} G_z \sum_{k_2, i_2} \sum_{k_3, i_3} \Rtens^{(\chan_z)}_{(k_1 i_1), (k_2 i_2), (k_3 i_3)} \, f_{(k_2 i_2), q_2^z} \, f_{(k_3 i_3), q_3^z},
    \label{eq:online_contraction}
\end{equation}
where $z$ streams through the sparse geometric COO list, fetching the Gaunt weight $G_z$, the test and interacting angular states $(q_1^z, q_2^z, q_3^z)$, and the physical channel $\chan_z$ \citep{hiemstra2026wigner}. Each test state $q_1$ collects the $\numG$ transitions routed to it. The polyatomic operator differs from the monatomic contraction in a single respect: the dense radial direction carries the combined radial--internal mode pair $(k, i)$. Each channel block of $\Rtens$ is therefore a dense cube of side $N_{KI} = (\kmax+1)(\imax+1)$ where the monatomic block has side $\kmax + 1$. The angular routing (COO list, channels, and slices) is untouched. The contraction kernels are parameterized by the length of the dense radial block and agnostic to its composition, so the monatomic kernels are reused verbatim with the flattened pair index $(k, i)$.

Because addition is commutative, the summation loops of \eqref{eq:online_contraction} can be rearranged. The three algorithms below evaluate the same sum and differ only in the order in which the sparse angular routing and the dense radial--internal sums are executed, and therefore in their memory profile. They are the orderings of the monatomic study (Section~4.5 of \citet{hiemstra2026wigner}) carried over with $\kmax + 1$ replaced by $N_{KI}$. Table~\ref{tab:contraction_orderings} collects their operation counts and intermediate footprints.

\subsubsection*{Naive streaming}

The sum \eqref{eq:online_contraction} is evaluated as written: the algorithm streams linearly through the sparse COO Gaunt list and executes the dense radial--internal triple sum for every geometric transition. The sparse routing is coupled directly to the dense arithmetic. The angular states $(q_2^z, q_3^z)$ jump pseudo-randomly across the spectral space as $z$ increments, so the coefficient fetches are cache-unfriendly.

\subsubsection*{Radial-first}

The radial--internal exchange is evaluated before any geometric routing. Contracting the dense physical tensor with the distribution coefficients for every trial pair $(q_2, q_3)$ and channel $\chan$ carried by the COO list builds the dense intermediate $\Psi$ of \citet{hiemstra2026wigner},
\begin{equation}
    \Psi^{(k_1 i_1)}_{q_2, q_3, \chan} = \sum_{k_2, i_2} \sum_{k_3, i_3} \Rtens^{(\chan)}_{(k_1 i_1), (k_2 i_2), (k_3 i_3)} \, f_{(k_2 i_2), q_2} \, f_{(k_3 i_3), q_3},
    \label{eq:contraction_radial}
\end{equation}
after which the operator follows in a single streaming pass over the COO list,
\begin{equation}
    Q_{(k_1 i_1), q_1} = \sum_{\{z \,\mid\, q_1^z = q_1\}} G_z \, \Psi^{(k_1 i_1)}_{q_2^z, q_3^z, \chan_z}.
    \label{eq:contraction_radial_stream}
\end{equation}
The nested radial--internal loops thereby leave the sparse execution path, at the cost of allocating $\Psi$.

\subsubsection*{Angular-first (sliced)}

The angular geometry is aggregated before the physical tensor is applied. For a fixed test state $q_1$ and channel $\chan$, the COO list is swept once to build the compact angular matrix $\Phi$ of \citet{hiemstra2026wigner},
\begin{equation}
    \Phi^{(\chan, q_1)}_{(k_2 i_2), (k_3 i_3)} = \sum_{\{z \,\mid\, \chan_z = \chan, \; q_1^z = q_1\}} G_z \, f_{(k_2 i_2), q_2^z} \, f_{(k_3 i_3), q_3^z},
    \label{eq:contraction_angular}
\end{equation}
which is then contracted against the dense channel block and summed over the $\numT$ physical channels,
\begin{equation}
    Q_{(k_1 i_1), q_1} = \sum_{\chan = 1}^{\numT} \sum_{k_2, i_2} \sum_{k_3, i_3} \Rtens^{(\chan)}_{(k_1 i_1), (k_2 i_2), (k_3 i_3)} \, \Phi^{(\chan, q_1)}_{(k_2 i_2), (k_3 i_3)}.
    \label{eq:contraction_angular_apply}
\end{equation}
The sparse lookups are decoupled from the dense radial--internal arithmetic, and the intermediate is a single $N_{KI} \times N_{KI}$ matrix.

\begin{table}[t]
\centering
\caption{Contraction orderings for the online evaluation of \eqref{eq:online_contraction}: $\numG$ is the Gaunt nonzero count, $\numS \leq \numG$ the number of unique $(\chan, q_1)$ slices, and $N_{KI} = (\kmax+1)(\imax+1)$ the side of the dense radial--internal block; the monatomic counts of \citet{hiemstra2026wigner} are recovered with $N_{KI} \to \kmax + 1$. All three factorized orderings retain the $\mathcal{O}(\lmax^5)$ geometric scaling of the factorization; the dense Cartesian baseline remains $\mathcal{O}(N^3) = \mathcal{O}(\kmax^3 \imax^3 \lmax^6)$.}
\label{tab:contraction_orderings}
\begin{tabular}{@{}llll@{}}
\toprule
\textbf{Ordering} & \textbf{Operations} & \textbf{Intermediate} & \textbf{Footprint} \\
\midrule
Naive streaming & $\mathcal{O}(\numG N_{KI}^3)$ & --- & --- \\
Radial-first & $\mathcal{O}(\numG N_{KI}^3) + \mathcal{O}(\numG N_{KI})$ & $\Psi$, all trial pairs & $\numG N_{KI}$ \\
Angular-first (sliced) & $\mathcal{O}(\numG N_{KI}^2) + \mathcal{O}(\numS N_{KI}^3)$ & $\Phi$, per $(\chan, q_1)$ slice & $N_{KI}^2$ \\
\bottomrule
\end{tabular}
\end{table}

The enlarged radial dimension does not overturn the balance between the orderings. The advantage of the angular-first ordering rests on the intermediate $\Phi$ residing in the highest cache levels. Its $N_{KI}^2$ footprint grows quadratically with the internal resolution but stays first-level-cache resident across accessible truncations, and leaving that regime would require $N_{KI} \approx 128$, an order of magnitude beyond the resolutions at which the polyatomic operator is used. Its margin over the radial-first ordering therefore widens with internal resolution rather than narrowing, as the measurements of Section~\ref{sec:results} confirm.

\section{Results}
\label{sec:results}

This section validates the discrete operator in three parts. First, we establish the spectral convergence of the tensorized quadrature (Section~\ref{sec:res_convergence}) and the recovery of the monatomic limit (Section~\ref{sec:res_mono}). Second, we verify the physical fidelity of the assembled operator: conservation of the collision invariants (Section~\ref{sec:res_conservation}), transient relaxation of the internal temperature (Section~\ref{sec:res_relax}), and the transport coefficients of N$_2$, CO, and H$_2$ against experiment (Section~\ref{sec:res_transport}). Finally, we benchmark storage compression and contraction performance (Section~\ref{sec:res_perf}).

\subsection{Spectral Convergence of the Polyatomic Quadrature}
\label{sec:res_convergence}

To verify the composition argument of Section~\ref{sec:quadfamilies} numerically, we measure the integration error of the physical tensor against an over-resolved reference. The tensor is assembled for the non-frozen DSMC kernel ($\inelprob = 1$, $\dofpoly = 2.01$) at the spectral truncation $\kmax = \lmax = \imax = 2$. The padding of the spatial coordinates is swept over $N_{\mathrm{pad}} \in \{0, 1, 2, 4, \ldots, 24\}$ against a reference tensor computed at padding $32$. The internal axes run at their exactness-bound node counts in both the sweep and the reference (Section~\ref{sec:quadfamilies}), so the internal-axis contribution is exact and the sweep isolates the spatial architecture. Conservation enforcement is disabled. Three interaction exponents are tested: the Maxwell ($\vhs = 0$) and hard-sphere ($\vhs = 1$) limits, and the N$_2$-fitted exponent ($\vhs = 0.533$). The global tensor error is the relative $\ell_\infty$ maximum norm
\begin{equation}
    \epsilon = \frac{\| \Rtens_{N_{\mathrm{pad}}} - \Rtens_{\mathrm{ref}} \|_{\ell_\infty}}{\| \Rtens_{\mathrm{ref}} \|_{\ell_\infty}}.
    \label{eq:conv_metric}
\end{equation}

\begin{figure}[t]
\centering
\includegraphics[width=0.62\textwidth]{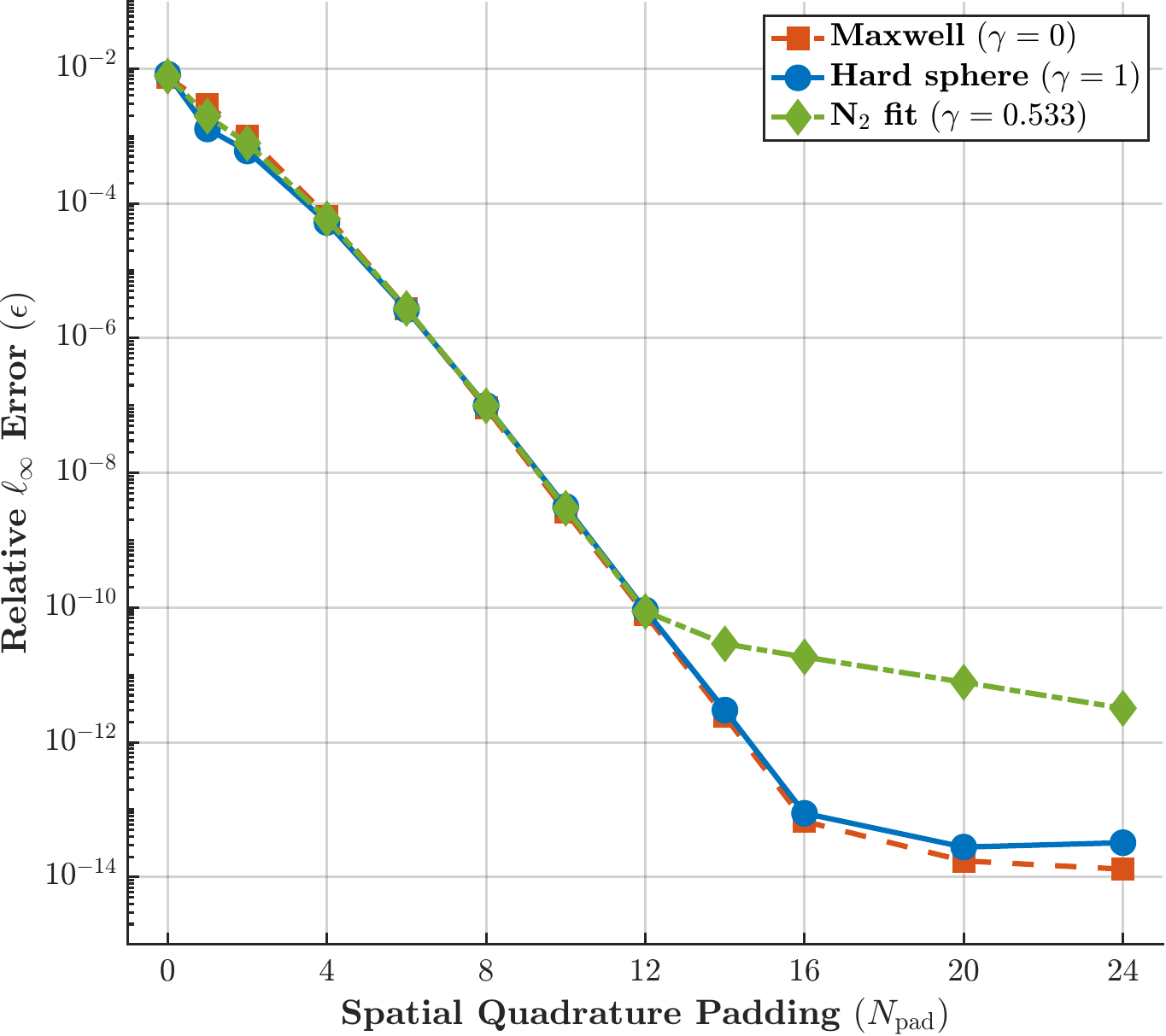}
\caption{Relative $\ell_\infty$ error \eqref{eq:conv_metric} of the physical tensor $\Rtens$ against the spatial quadrature padding, for the non-frozen DSMC kernel ($\inelprob = 1$, $\dofpoly = 2.01$) at $\kmax = \lmax = \imax = 2$, with reference padding $32$ and the internal axes at their exactness-bound node counts in sweep and reference. The integer-exponent limits reach machine precision, while the fractional N$_2$ exponent leaves a slow algebraic endpoint tail well below working accuracy.}
\label{fig:quadrature_convergence}
\end{figure}

As shown in Figure~\ref{fig:quadrature_convergence}, all three kernels converge monotonically and geometrically at a common rate toward the high-resolution reference, with the two integer-exponent limits reaching machine plateaus below $10^{-13}$ at $N_{\mathrm{pad}} = 16$. The 10--20 point padding experience of the monatomic architecture thus carries over unchanged \citep{hiemstra2026wigner}. The fractional N$_2$ exponent follows the same geometric decay before crossing into a slow algebraic endpoint tail, four to five orders below working accuracy, originating from the fractional power $\urel^{\vhs}$ left after the Duffy maps. The common decay rate shows that the Duffy patch pair regularizes the cone singularity uniformly in the interaction exponent.

\subsection{Validation in the Monatomic Limit}
\label{sec:res_mono}

To verify that the nine-dimensional polyatomic machinery introduces no spurious cross-coupling, we validate the monatomic limit in two complementary ways: exactly, at the truncation $\imax = 0$ where the Dirac-collapse rule of Section~\ref{sec:quadfamilies} is active, and continuously, along a sweep of the internal dimension $\dofpoly \to 0$.

\subsubsection{Exact Reproduction at the Collapse Point}
For the non-frozen Maxwell kernel ($\vhs = 0$, $\inelprob = 1$, $\dofpoly = 2$), the collapsed operator agrees with an independent monatomic Maxwell build to a relative $\ell_\infty$ difference of $2.2 \times 10^{-15}$, with the equilibrium collision frequency matching the monatomic VHS normalization identically. The full 27-mode Wang Chang--Uhlenbeck (WCU) eigenvalue spectrum \citep{wangchang:1952,kumar:1966} ($\kmax = \lmax = 2$) matches the analytic values to $4 \times 10^{-14}$ or better, with the five collision invariants exactly zero. Conservation is native at the collapse point: the raw invariant residuals are $4.1 \times 10^{-15}$ before the projection of Section~\ref{sec:conservation} is applied.

\subsubsection{Continuous Approach in the Internal Dimension}
\begin{figure}[t]
\centering
\includegraphics[width=\textwidth]{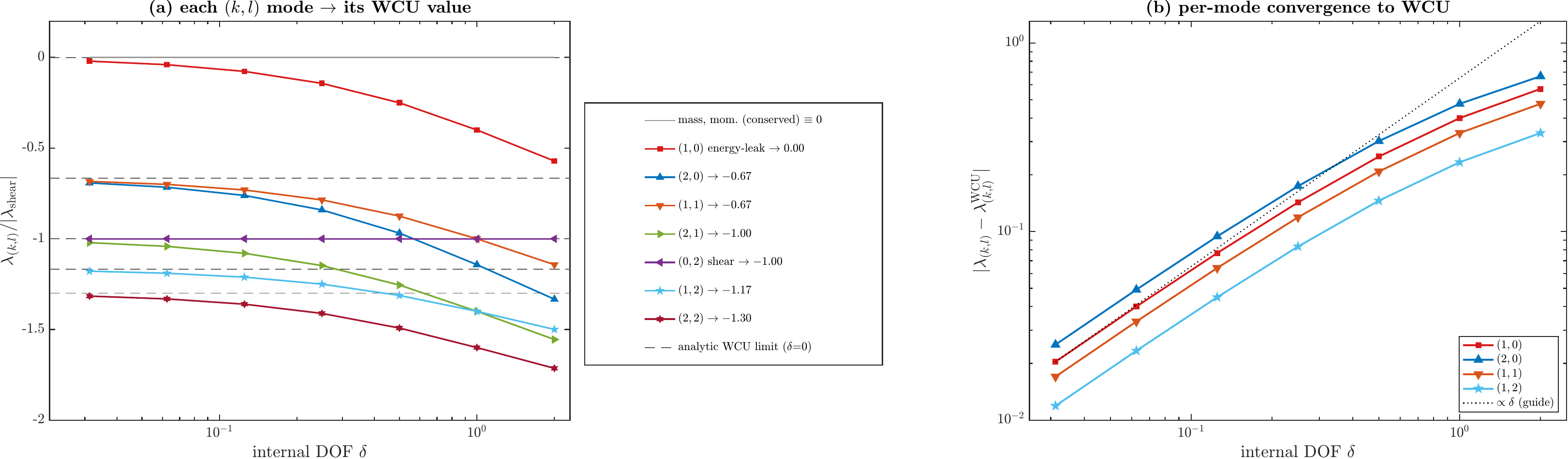}
\caption{Convergence of the polyatomic eigenvalue spectrum to the Wang Chang--Uhlenbeck (monatomic) limit as the internal degrees of freedom vanish, for the non-frozen Maxwell kernel ($\vhs = 0$, $\inelprob = 1$) at $\kmax = \lmax = \imax = 2$ with $(16,16)$ spatial padding, internal axes at their exactness-bound node counts, and $\dofpoly \in \{2, 1, \ldots, 1/32\}$. The plotted values are the Rayleigh--Ritz eigenvalues of the $i = 0$, $m = 0$ principal sub-block of the linearized operator, taken separately for each $l$. Left: the shear-normalized values $\lambda_{(k,l)}/|\lambda_{\mathrm{shear}}|$ approach the analytic WCU values (dashed) monotonically; mass and momentum are exact at every $\dofpoly$, while the shear mode sits at $-1$ by construction of the normalization. Right: the per-mode deviation decays proportionally to $\dofpoly$. The $\dofpoly = 0$ endpoint of the sweep is singular ($\polyalpha = -1$) and is anchored by the analytic values.}
\label{fig:wcu_limit}
\end{figure}

The operator is assembled for the non-frozen Maxwell kernel ($\vhs = 0$, $\inelprob = 1$) at $\kmax = \lmax = \imax = 2$ with spatial padding $(16, 16)$ on the radial and angular axes and the four internal axes at their exactness-bound node counts of Table~\ref{tab:quadrature_bounds}. At this truncation these are $N_{\Ipre} = N_{\Jpre} = N_{\partkin} = 4$ and $N_{\partint} = 3$, the convention of every production run below. The internal dimension is swept geometrically over $\dofpoly \in \{2, 1, \tfrac{1}{2}, \ldots, \tfrac{1}{32}\}$, approaching the singular endpoint $\polyalpha = \dofpoly/2 - 1 \to -1$ at which the internal equilibrium weight collapses onto a Dirac measure. The conservation enforcement of Section~\ref{sec:conservation} pins the five invariant rows exactly at every $\dofpoly$. The quantities plotted in Figure~\ref{fig:wcu_limit} are the eigenvalues of the $i = 0$, $m = 0$ principal sub-block of the linearized operator, taken separately for each $l$: that subspace is not invariant under the polyatomic operator, so they are Rayleigh--Ritz values rather than eigenvalues of the discrete operator itself. Every dynamic mode of the shear-normalized set converges monotonically to its WCU value, with a per-mode deviation decaying proportionally to $\dofpoly$.

Because the spectrum is normalized by $|\lambda_{\mathrm{shear}}|$, the shear mode $(0,2)$ sits at $-1$ by construction, and its accuracy is established separately in absolute terms. We measure $\nu_0$ from a loss-only build of the same operator, obtained by zeroing the gain kernels, so it involves no eigenvalue of the operator under test. At $\vhs = 0$ the linearized loss operator is $-\nu_0 I$ on every non-density mode. Against it the shear eigenvalue is $\lambda_{\mathrm{shear}} = -\nu_0/2$ at every $\dofpoly$, matching the analytic value \eqref{eq:lambda_shear} of Appendix~\ref{app:maxwell} to below $10^{-13}$, and $\nu_0$ itself matches its analytic value to the same order. The mode is an exact eigenvector of the full discrete operator and does not couple to the internal sector, as the derivation of Appendix~\ref{app:maxwell} assumes. Over the same sweep the energy-exchange eigenvalue satisfies $\lambda_{\mathrm{exch}} = -\nu_0 (3 + \dofpoly)/(3 + 2\dofpoly)$ \eqref{eq:lambda_bulk}, again to below $10^{-13}$, and its ratio to the shear eigenvalue matches the analytic $\lambda_{\mathrm{exch}}/\lambda_{\mathrm{shear}} = 10/7$ of Appendix~\ref{app:maxwell} at $\dofpoly = 2$. This eigenvalue carries internal-energy content, so it lies outside the $i = 0$ sub-block of Figure~\ref{fig:wcu_limit}. The complementary infinite-order hard-sphere viscosity benchmark is deferred to the transport study of Section~\ref{sec:res_transport}.

\subsection{Conservation Laws and the Polyatomic Null-Space}
\label{sec:res_conservation}

The conservation embedding of Section~\ref{sec:conservation} must hold under nonlinear evolution, where invariant leakage compounds across composed operator evaluations. The static enforcement itself is verified separately (twelve parameterized configurations, all exact). Here we track the invariants over the full transient of the two-temperature relaxation experiment of Section~\ref{sec:res_relax}: nonlinear evolution from $(T_v, T_I) = (4/3, 1)$ under the DSMC kernel ($\dofpoly = 2$, $\inelprob = 1$, $\vhs = 0.5$), integrated by a second-order Runge--Kutta (Heun) scheme. A Maxwell control ($\vhs = 0$), quadrature-exact by construction, is run alongside on the same grid and integrator to isolate the time stepper's contribution to the residual drift.

\begin{figure}[t]
\centering
\includegraphics[width=0.62\textwidth]{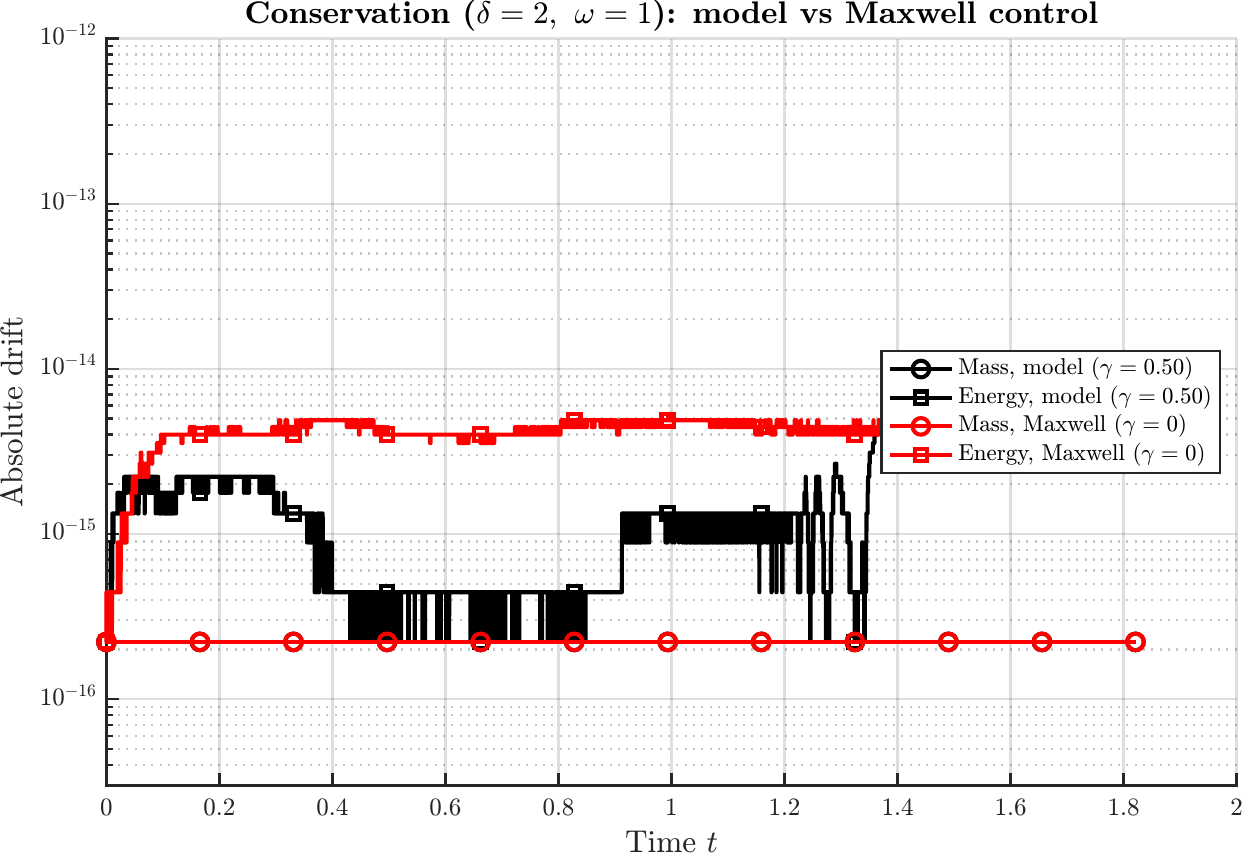}
\caption{Invariant drift over the two-temperature relaxation transient of Section~\ref{sec:res_relax}: absolute deviation of mass and total energy from their initial values for the production DSMC model ($\vhs = 0.5$, black) and a Maxwell control ($\vhs = 0$, red) on the same grid and second-order Runge--Kutta integrator, at $\dofpoly = 2$, $\inelprob = 1$, $\kmax = \lmax = \imax = 2$, with $(10,10)$ spatial padding and internal axes at their exactness-bound node counts. The mass drift is exactly zero and is plotted at the $10^{-16}$ floor. The energy drift of both operators stays below $3 \times 10^{-14}$.}
\label{fig:polyrelax_conservation}
\end{figure}

\textbf{Measured drift.} The mass drift is exactly zero for both operators throughout the transient, and the total-energy drift stays below $3 \times 10^{-14}$ (Figure~\ref{fig:polyrelax_conservation}). We quote a bound rather than a value: the offline build reduces in nondeterministic thread order, so the trailing digits of quantities at this level are not reproducible run to run.

\textbf{Attribution.} The comparable magnitudes attribute the residual to the time stepper rather than the quadrature: the null-space embedding of Section~\ref{sec:conservation} holds exactly, independently of both.

\textbf{Unenforced residual.} Because the embedding zeroes the mass and momentum rows outright, their exact conservation is a property of the construction. The quadrature attains it independently: with enforcement disabled, the mass and momentum rows of the production operators of Section~\ref{sec:res_transport} sit at assembly roundoff, below $2 \times 10^{-15}$ relative, and the energy combination that the projection removes is below $10^{-13}$ relative. The embedding therefore corrects a residual that is already at the level of the arithmetic, rather than masking a quadrature error.

\textbf{Equilibrium residual.} The discrete operator also annihilates the reference equilibrium. Across the production configurations the relative residual $\|\Ctens_{\state, 0, 0}\|$ lies between $10^{-11}$ and $10^{-8}$, with the base kernels at the lower end and the modulated extended kernel at the upper end. This is the discrete counterpart of detailed balance, and it is not implied by the invariant projection: the projection constrains test rows, whereas this residual is a statement about the trial state.

\subsection{Transient Relaxation of Internal Temperatures}
\label{sec:res_relax}

To validate the inelastic exchange physics of the partition integrals, we evolve a spatially homogeneous gas from thermal non-equilibrium and track the translational and internal temperatures $T_v$ and $T_I$ against the analytic relaxation laws. The gas is initialized as a product of Maxwellians at $(T_{v,0}, T_{I,0}) = (4/3, 1)$ and evolved nonlinearly at $\kmax = \lmax = \imax = 2$ with $(10, 10)$ spatial padding and internal axes at their exactness-bound node counts, the operator rescaled to unit equilibrium collision frequency $\nu_0 = 1$.

\textbf{Internal-dimension sweep.} For the Maxwell kernel ($\vhs = 0$, $\inelprob = 1$), the internal dimension $\dofpoly$ fixes both the equilibrium temperature and the exchange rate through the equipartition law and the Landau--Teller relaxation law \citep{landau:1936}
\begin{equation}
    T_{\eq} = \frac{3 T_{v,0} + \dofpoly \, T_{I,0}}{3 + \dofpoly}, \qquad
    \frac{d}{dt}(T_v - T_I) = -\lambda_{\mathrm{LT}} \, (T_v - T_I), \qquad
    \lambda_{\mathrm{LT}} = \nu_0 \, \frac{3 + \dofpoly}{3 + 2\dofpoly},
    \label{eq:landau_teller}
\end{equation}
the latter exact for Maxwell molecules, whose exchange rate is speed-independent. \citet{landau:1936} establish the exponential relaxation law itself. The rate $\lambda_{\mathrm{LT}}$ in terms of the internal degrees of freedom is specific to the Borgnakke--Larsen kernel and is derived in Appendix~\ref{app:maxwell}.

As shown in Figure~\ref{fig:polyrelax_di}, the runs $\dofpoly \in \{2, 3, 4\}$ relax to the rational equilibria $T_{\eq} = 6/5$, $7/6$, $8/7$, and the temperature gap decays exponentially at the Landau--Teller rates $\lambda_{\mathrm{LT}} = 5/7$, $2/3$, $7/11$. The energy-exchange eigenvalue of the linearized operator matches these rates to relative errors below $10^{-13}$. The plotted transients are themselves single exponentials: a log-linear fit recovers the same rates to $9.6 \times 10^{-7}$, window-independent over five decades of amplitude, the residual being exactly the second-order amplification error $(\lambda \Delta t)^2/6$ of the Heun integrator.

\textbf{Inelastic-probability sweep.} Sweeping instead the inelastic probability $\inelprob \in \{0.25, 0.5, 0.75, 1\}$ for the production kernel ($\vhs = 0.5$, $\dofpoly = 2$), all four runs relax to the same equilibrium $T_{\eq} = 6/5$, since energy conservation is independent of the channel split. The exchange eigenvalue is strictly proportional to $\inelprob$, its ratio to $\inelprob$ constant at $26.357$ across the sweep to roundoff (Figure~\ref{fig:polyrelax_omega}). The proportionality is structural: the translational- and internal-energy rows of the frozen operator vanish separately, so the two-mode exchange block is exactly $\inelprob$ times its value at $\inelprob = 1$. Fitted instead to the plotted transients the ratio is $5.2\%$ higher and drifts by $0.55\%$ across the sweep, because at $\vhs = 0.5$ the exchange block is no longer an invariant subspace and the observable decays as a mixture of the neighboring $l = 0$ modes. The late-time plateau of $|T_v - T_I|$, near $2.5 \times 10^{-4}$ across the $\inelprob \geq 0.5$ runs, is the accuracy limit of the finite spectral truncation. It is unchanged when the spatial padding is refined to $(16,16)$ or the time step is halved, but falls by roughly an order of magnitude, to $2.7 \times 10^{-5}$, when the basis is enriched to $\kmax = \imax = 3$. It lies outside the early-time window of the rate fits.

\begin{figure}[t]
\centering
\includegraphics[width=\textwidth]{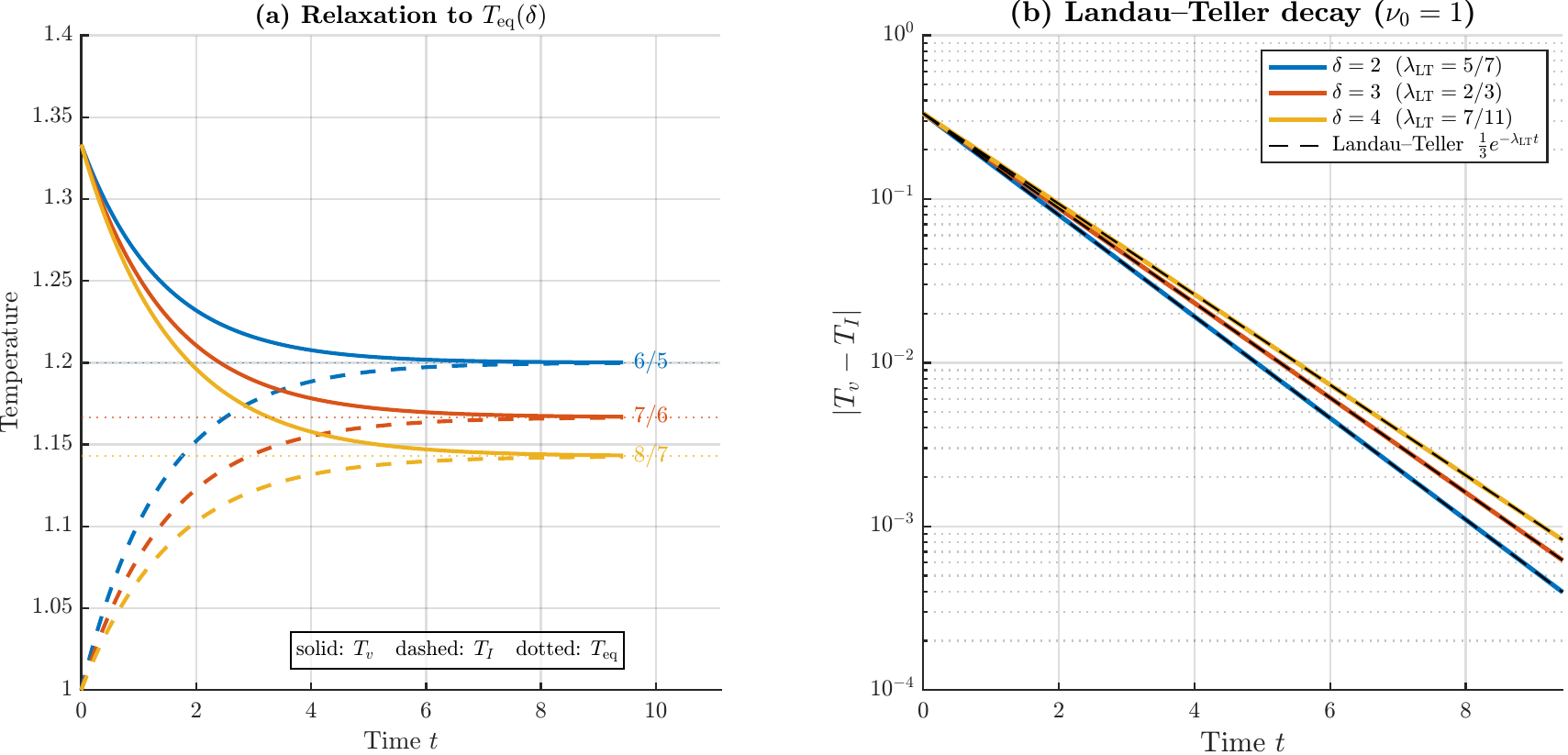}
\caption{Two-temperature relaxation under the Borgnakke--Larsen Maxwell kernel ($\vhs = 0$, $\inelprob = 1$) for internal dimension $\dofpoly \in \{2, 3, 4\}$, at $\kmax = \lmax = \imax = 2$ with $(10,10)$ spatial padding and internal axes at their exactness-bound node counts, from the product-Maxwellian initial state $(T_{v,0}, T_{I,0}) = (4/3, 1)$, with unit equilibrium collision frequency $\nu_0 = 1$. Left: translational (solid) and internal (dashed) temperatures relax to the rational equilibria $T_{\eq} = (3 T_{v,0} + \dofpoly T_{I,0})/(3 + \dofpoly)$ (dotted). Right: the temperature gap $|T_v - T_I|$ decays exponentially at the Landau--Teller rate \eqref{eq:landau_teller} (dashed overlay).}
\label{fig:polyrelax_di}
\end{figure}

\begin{figure}[t]
\centering
\includegraphics[width=\textwidth]{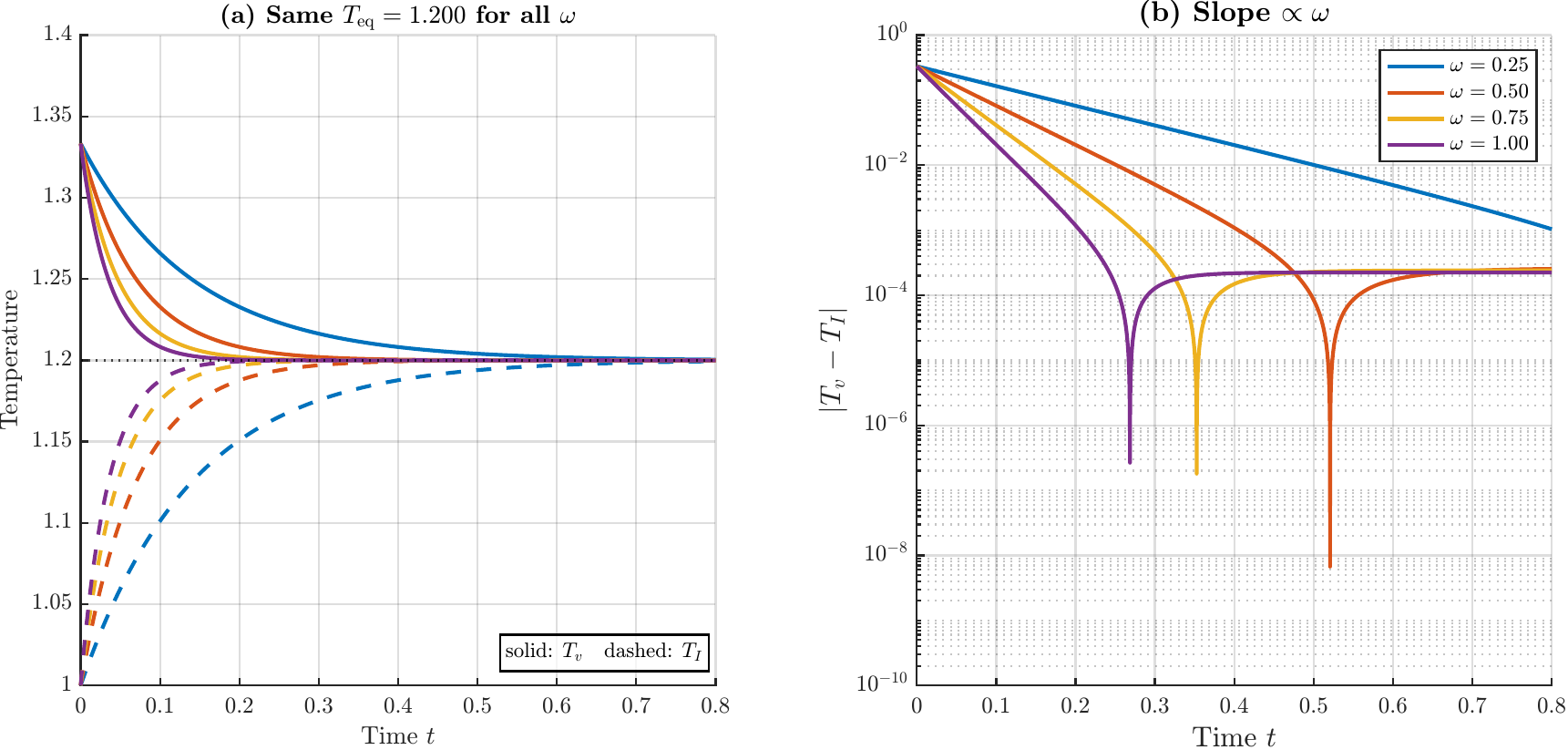}
\caption{Two-temperature relaxation under the production DSMC kernel ($\vhs = 0.5$, $\dofpoly = 2$) for inelastic probability $\inelprob \in \{0.25, 0.5, 0.75, 1\}$, with resolution, padding, and initial state as in Figure~\ref{fig:polyrelax_di}; the frozen and non-frozen channels are built once and convex-blended in $\inelprob$. Left: all four runs reach the same equilibrium $T_{\eq} = 6/5$. Right: the exchange eigenvalue is strictly proportional to $\inelprob$. The late-time plateaus of the $\inelprob \geq 0.5$ runs, all within $[2.2, 2.6] \times 10^{-4}$, are the truncation-limited discrete equilibrium discussed in the text. The downward spikes mark the change of sign of $T_v - T_I$.}
\label{fig:polyrelax_omega}
\end{figure}

\subsection{Polyatomic Transport Coefficients}
\label{sec:res_transport}

Bulk viscosity is the defining macroscopic signature of a polyatomic gas: it vanishes identically in the monatomic limit. Together with the Prandtl number, it is the quantity against which the extended DSMC kernel \eqref{eq:dsmc_kernel} was calibrated by \citet{djordjic:2023}. We therefore assemble the operator for three real gases and compare its transport coefficients against experiment: first at the published kernel parameters, with nothing fitted, and then under a calibration of the kernel's free parameter pair within our own operator.

\subsubsection{Coefficient Extraction and Setup}
\label{sec:res_transport_setup}
The coefficients are the first-order production coefficients of the operator linearized about absolute equilibrium: the shear coefficient $P_\sigma$ as the diagonal $l = 2$ element of the $(k,i) = (0,0)$ stress mode, the heat-flux pair as the $l = 1$ block on the translational and internal flux modes, and the bulk coefficient $P_\Pi$ as the non-zero eigenvalue of the $l = 0$ energy-exchange block, the block that is empty in the monatomic operator. These feed the $17$-moment relations of \citet{djordjic:2023}, Eqs.~(39) and (58). The $l = 1$ block is expressed in their moments by $q = d_q c_q$ and $s = d_s c_s$. Projecting $\tfrac{1}{2}|\mathbf{c}|^2 c_i$ and $\Ipre c_i$ onto the orthonormal modes gives $d_q/d_s = \sqrt{5/\dofpoly}$, the $l = 1$ counterpart of the ratio $\sqrt{3/\dofpoly}$ that fixes the $l = 0$ energy invariant. This is the closure through which their Table-4 parameters were calibrated. Its sensitivity to the Chapman--Enskog order is measured in Section~\ref{sec:res_transport_order}. The two reported quantities are the Prandtl number
\begin{equation}
    \prandtl = \frac{c_p \, \shearvisc}{\heatcond}, \qquad c_p = \frac{5 + \dofpoly}{2},
    \label{eq:prandtl}
\end{equation}
with $c_p$ in units of the specific gas constant, and the bulk-to-shear viscosity ratio $\bulkvisc/\shearvisc$. We write $\bulkvisc$ for the bulk viscosity that they denote $\nu$, since $\polyalpha$ is the internal-Laguerre parameter throughout this paper. For each gas, $\dofpoly$ and the VHS exponent $\vhs$ are those of their Table~1. All operators are assembled at $\kmax = \lmax = \imax = 2$ with $(16, 16)$ spatial padding on the radial and angular axes, the internal axes at their exactness-bound node counts, $\numaux = 24$ auxiliary Laplace nodes, and conservation enforcement active. The spatial padding is the value at which the integer-exponent limits of Section~\ref{sec:res_convergence} plateau and the fractional-exponent tail of these gases is below $2 \times 10^{-11}$. Internal padding is immaterial: internal paddings $0$ and $8$ agree to within $2 \times 10^{-14}$ on all channel branches.

\subsubsection{Transport Coefficients at the Published Parameters}
\label{sec:res_transport_pred}
Assembled at the published parameter values $(\inelprob, \modamp, \modexp, \modampf, \modexpf)$ of \citet{djordjic:2023}, Table~4, the operator reproduces both measured coefficients for every gas to better than $0.0002\%$. The Prandtl numbers come out at $0.71700002$ against $0.717$ for N$_2$, $0.74300002$ against $0.743$ for CO, and $0.68600002$ against $0.686$ for H$_2$. The bulk-to-shear ratios come out at $0.72999897$ against $0.730$, $0.54999981$ against $0.550$, and $30.00005279$ against $30.0$.

The agreement with experiment is inherited, not predicted. The Table-4 parameters of \citet{djordjic:2023} were themselves obtained by solving $\prandtl = \prandtl_{\mathrm{exp}}$ and $\bulkvisc/\shearvisc = (\bulkvisc/\shearvisc)_{\mathrm{exp}}$ for $(\inelprob, \modampf)$ through the same first-order closure we use in Section~\ref{sec:res_transport_setup}. Evaluating at their parameters and recovering the experimental values therefore tests our quadrature against their calibration, not the kernel against experiment. What the comparison establishes independently is a code-to-code verification: a nine-dimensional quadrature of the kernel reproduces their symbolically evaluated production coefficients to within $0.0002\%$. That margin is four orders of magnitude below the systematic carried by the extraction itself, since resolving the Chapman--Enskog order moves the bulk-to-shear ratio by several percent (Section~\ref{sec:res_transport_order}).

One comparison in this section is parameter-free. In the pure frozen limit $\inelprob = 0$ the operator reproduces the closed-form frozen Prandtl number of \citet{djordjic:2023}, their Eq.~(42), to $2 \times 10^{-7}\%$ for all three gases, with the residual set by the spatial padding rather than by a floor. The frozen channel is therefore verified against an analytic target independently of the non-frozen one, and the agreement at the published parameters above holds across the whole range of the inelastic fraction, from $\inelprob = 0.010$ for H$_2$ to $0.541$ for CO.

\subsubsection{Calibration and Recovery of the Published Parameters}
\label{sec:res_transport_fit}
\begin{figure}[t]
\centering
\includegraphics[width=\textwidth]{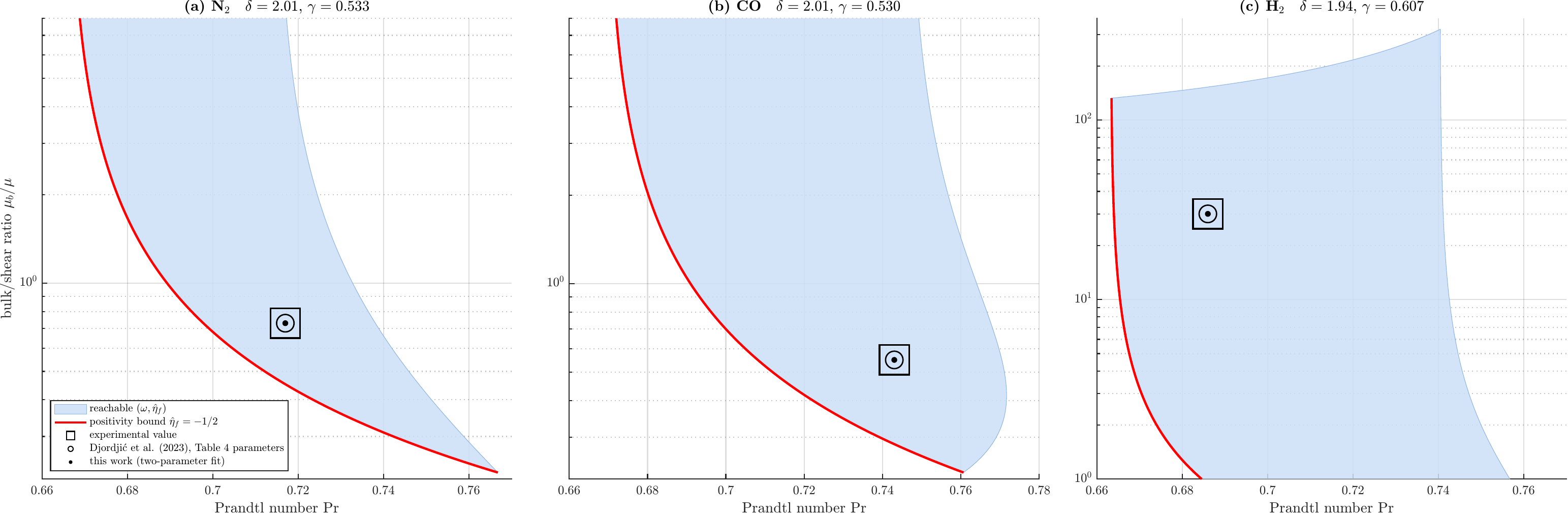}
\caption{Reachability of the experimental transport targets in the $(\prandtl, \bulkvisc/\shearvisc)$ plane for the extended DSMC kernel \eqref{eq:dsmc_kernel}, at $\kmax = \lmax = \imax = 2$, with $(16,16)$ spatial padding, internal axes at their exactness-bound node counts, and $\numaux = 24$ auxiliary Laplace nodes. The shaded region is the image of the two-parameter family $(\inelprob, \modampf)$, sampled on a $420 \times 420$ grid over $\inelprob \in [2 \times 10^{-3}, 1]$ and $\modampf \in [-1/2, 3]$, with the remaining modulation parameters fixed at Table~4 of \citet{djordjic:2023}; the upper edge of the shaded region in $\prandtl$ is therefore the sampling limit rather than a bound of the model. The red curve is the positivity boundary $\modampf = -1/2$. Each panel carries three coincident markers, drawn back to front: the experimental value (large open square), the operator at the published parameters of \citet{djordjic:2023}, Table~4 (open circle), and the calibration of this work (filled dot), which overlap in every panel.}
\label{fig:transport_fits}
\end{figure}

We now free the published parameters, calibrating the kernel within our own operator: holding $(\modamp, \modexp, \modexpf)$ at their published values, the pair $(\inelprob, \modampf)$ is fitted to the experimental Prandtl number and bulk-to-shear ratio of each gas. Both targets are reproduced exactly for all three gases, at $(\inelprob, \modampf) = (0.3121, -0.2078)$ for N$_2$, $(0.5405, +0.5701)$ for CO, and $(0.0101, -0.1339)$ for H$_2$ (Table~\ref{tab:transport}). The fitted pair enters only through the convex blend of the two channel operators and the amplitude linearity of Section~\ref{sec:laplace}, so the fit is a two-parameter solve on pre-built tensors with nothing reassembled. Each fit is a Newton solve on the two-dimensional residual with a finite-difference Jacobian, and all three respect the kernel positivity bound $\modampf \geq -1/2$ imposed on \eqref{eq:dsmc_kernel}. Figure~\ref{fig:transport_fits} places every experimental target inside the admissible image of the two-parameter family, away from its boundary. The reachable Prandtl numbers span $[0.666, 0.767]$ for N$_2$, $[0.669, 0.772]$ for CO, and $[0.663, 0.764]$ for H$_2$. The reachable bulk-to-shear ratio of H$_2$ spans more than three orders of magnitude, so its experimental ratio of $30$ is met well inside the region.

The fitted probabilities are physically sensible. With $\xi_{\mathrm{tr}} = 3 + \vhs$, Eq.~(44) of \citet{djordjic:2023} relates the inelastic probability to the internal collision number by $Z_{\mathrm{int}} = (3 + \vhs + \dofpoly)/[(3 + \vhs)\inelprob]$, and our fitted values give $Z_{\mathrm{int}} = 5.03$ for N$_2$, $2.90$ for CO, and $152$ for H$_2$. The N$_2$ value sits on the standard room-temperature rotational collision number $Z_{\mathrm{rot}} \approx 5$, and H$_2$'s large value reflects its wide rotational spacing.

The calibration also recovers the published parameters. Fitted independently through our own quadrature, the pair $(\inelprob, \modampf)$ reproduces the values of \citet{djordjic:2023}, Table~4, to within $4 \times 10^{-6}$ in relative terms for all three gases, a like-for-like check of two independent routes to the same calibration.

\subsubsection{Sensitivity to the Chapman--Enskog Order}
\label{sec:res_transport_order}
The first-order extraction of Section~\ref{sec:res_transport_setup} reads three fixed lead modes. To measure what it omits, we replace it by a full inversion of each angular block, Schur-reduced onto the same macroscopic moments, and resolve the radial (Sonine) truncation to $\kmax = 3$, at which the successive differences are below $10^{-3}$ in relative terms for both $\prandtl$ and $\bulkvisc/\shearvisc$.

The two quantities behave differently. The Prandtl number is robust: over all three gases and both parameter sets the resolved value differs from the first-order one by at most $0.26\%$. The bulk-to-shear ratio is not: the resolved value is higher by $4.5\%$ for H$_2$, $7.9\%$ for N$_2$, and $13.3\%$ for CO. The asymmetry is structural. The resolved bulk coefficient is the slowest non-zero eigenvalue of an $l = 0$ block that grows from two to eleven modes, and by eigenvalue interlacing that coefficient can only decrease as modes are added, so the ratio can only increase.

Re-fitting $(\inelprob, \modampf)$ under the resolved extraction therefore moves the parameters: $(0.324597, -0.254261)$ for N$_2$, $(0.603536, +0.512195)$ for CO, and $(0.010687, -0.105564)$ for H$_2$. All three still reproduce both experimental targets exactly and still satisfy the positivity bound $\modampf \geq -1/2$, so the reachability conclusion of Section~\ref{sec:res_transport_fit} is independent of the extraction; the fitted values are not. We report the first-order fits because they are calibrated through the same closure as the parameters they are compared against.

The resolved bulk coefficient is the spectral gap of an $l = 0$ block that grows with the internal truncation, so it depends on $\imax$; the first-order coefficient reads only $i = 0, 1$ and does not. Resolving $\imax$ to $4$ at $\kmax = 3$ leaves the ratio within $2 \times 10^{-3}$ in relative terms of its value at the production truncation $\imax = 2$. The successive differences contract geometrically for H$_2$, at ratio $0.27$, but stall near $10^{-4}$ for N$_2$ and CO, so the residual drift is bounded rather than resolved to a limit. In every case it is an order of magnitude or more below the first-order-to-resolved gap, so that comparison is not at risk from the internal truncation.

\subsubsection{Robustness and Temperature Scope}
Both comparisons are insensitive to the two numerical choices they depend on. The legacy algebraic path evaluates the fractional modulation pointwise at the tensorized nodes, without the auxiliary representation of Section~\ref{sec:laplace} and at its shipped internal padding of $8$. Recomputing all quantities on it changes $\prandtl$ by at most $1.6 \times 10^{-5}$ and $\bulkvisc/\shearvisc$ by at most $1.0 \times 10^{-3}$ in relative terms across the three gases. These differences are dominated by the legacy path's own residual error. At these truncations the auxiliary Laplace construction is a convergence guarantee rather than a large correction. Varying the truncation over $(\kmax, \imax) \in \{(1,2), (2,2), (3,2), (2,1), (2,3)\}$ on the production path changes both $\prandtl$ and $\bulkvisc/\shearvisc$ by at most $2 \times 10^{-9}$ in relative terms. The extraction reads fixed lead modes and is therefore independent of the truncation by construction, so this sweep measures the quadrature rather than the Chapman--Enskog order (Section~\ref{sec:res_transport_order}). The residual $\kmax$ dependence is spatial: the radial rules integrate $\urel^{\vhs}$, which is not polynomial, so their node counts never drop out entirely. Along the internal axis the three truncations agree to roundoff, because the internal rules are exact at every $\imax$ (Section~\ref{sec:quadfamilies}).

The temperature dependence of the transport coefficients is carried entirely by the variable-hard-sphere power law: $\vhs$ is fixed by the measured viscosity law $\shearvisc \propto T^{s_{\mathrm{visc}}}$ \citep{djordjic:2021}, and the modulation factors are dimensionless energy fractions. The ratios of Table~\ref{tab:transport} are therefore temperature-invariant within the calorically perfect regime, and a single fitted pair $(\inelprob, \modampf)$ holds across it. A measured drift of $\prandtl$ or $\bulkvisc/\shearvisc$ with temperature marks the boundary of the kernel family rather than a failure of the fit.

\subsubsection{Monatomic Viscosity Benchmark}
A full inversion of the $l = 2$ block, the resolved extraction of Section~\ref{sec:res_transport_order}, reproduces the classical Chapman--Enskog viscosity factors in the monatomic limit, the benchmark deferred from Section~\ref{sec:res_mono}. We report the shear viscosity correction factor $f_\mu(\kmax)$, the ratio of the shear viscosity obtained by inverting the full $l = 2$ block at truncation $\kmax$ to its first-order single-mode value. It is computed at the internal truncation $\imax = 0$, where the Dirac-collapse rule of Section~\ref{sec:quadfamilies} pins $\partkin = 1$ and the frozen channel is the monatomic elastic process exactly. For Maxwell molecules the factor is unity at every truncation, $\max |f_\mu - 1| = 0$, as the Sonine modes are eigenfunctions of the linearized Maxwell operator. For frozen hard spheres it converges monotonically to the Pekeris value $1.016034$ \citep{pekeris:1957}. The first correction reproduces the classical second-approximation ratio $205/202$ to every digit printed, and the successive Cauchy differences fall by roughly one order per level to a residual of $2.8 \times 10^{-5}$ at $\kmax = 3$ (Table~\ref{tab:viscosity_factor}). The infinite-order value is therefore a converged truncation rather than an extrapolation. Both cases are base kernels, $\modamp = \modampf = 0$, for which the auxiliary correction vanishes identically and the two quadrature paths coincide.

\begin{table}[t]
\centering
\caption{Infinite-order Chapman--Enskog shear viscosity correction factor $f_\mu(\kmax)$ for frozen hard spheres ($\vhs = 1$, $\inelprob = 0$) at $\dofpoly = 2$, $\lmax = 2$, $\imax = 0$, and $(16,16)$ spatial padding. The $\kmax = 1$ entry is the exact second-approximation ratio $205/202 = 1.014851485148515$, and the sequence converges to the Pekeris infinite-order value $1.016034$.}
\label{tab:viscosity_factor}
\begin{tabular}{@{}ccc@{}}
\toprule
$\kmax$ & $f_\mu$ & Cauchy difference \\
\midrule
$0$ & $1.000000000000000$ & --- \\
$1$ & $1.014851485148515$ & $1.485 \times 10^{-2}$ \\
$2$ & $1.015878911770992$ & $1.027 \times 10^{-3}$ \\
$3$ & $1.016005918782137$ & $1.270 \times 10^{-4}$ \\
\bottomrule
\end{tabular}
\end{table}

\begin{table}[t]
\centering
\caption{Polyatomic transport coefficients for the extended DSMC kernel \eqref{eq:dsmc_kernel} at $\kmax = \lmax = \imax = 2$, with $(16,16)$ spatial padding, internal axes at their exactness-bound node counts, $\numaux = 24$ auxiliary Laplace nodes, conservation enforcement active, and the internal degrees of freedom $\dofpoly$ and the VHS exponent $\vhs$ taken per gas from Table~1 of \citet{djordjic:2023}. Computed columns: the Prandtl number \eqref{eq:prandtl} and the bulk-to-shear viscosity ratio, extracted from the first-order production coefficients of the linearized operator (Section~\ref{sec:res_transport_setup}) at the fitted pair $(\inelprob, \modampf)$. The pair is calibrated within our operator against the experimental values of their Tables~2 and 3, at fixed $(\modamp, \modexp, \modexpf)$ with $\modexpf = 0.300$ for N$_2$ and $0.965$ for CO and H$_2$. The fit is a two-parameter solve on two targets, so the residual is zero by construction and the table reports a calibration rather than an independent test, with the pair quoted to four decimals since the first-order extraction carries a systematic of several percent on $\bulkvisc/\shearvisc$ (Section~\ref{sec:res_transport_order}).}
\label{tab:transport}
\resizebox{\textwidth}{!}{
\begin{tabular}{@{}lcccccccccc@{}}
\toprule
& & & \multicolumn{3}{c}{Prandtl number $\prandtl$} & \multicolumn{3}{c}{Bulk-to-shear ratio $\bulkvisc/\shearvisc$} & \multicolumn{2}{c}{Fitted parameters} \\
\cmidrule(lr){4-6} \cmidrule(lr){7-9} \cmidrule(lr){10-11}
Gas & $\dofpoly$ & $\vhs$ & Experiment & Computed & Deviation & Experiment & Computed & Deviation & $(\inelprob, \modampf)$ & Positivity \\
\midrule
N$_2$ & $2.01$ & $0.533$ & $0.717$ & $0.717000$ & $0.00\%$ & $0.730$ & $0.730000$ & $0.00\%$ & $(0.3121,\ -0.2078)$ & satisfied \\
CO    & $2.01$ & $0.530$ & $0.743$ & $0.743000$ & $0.00\%$ & $0.550$ & $0.550000$ & $0.00\%$ & $(0.5405,\ +0.5701)$ & satisfied \\
H$_2$ & $1.94$ & $0.607$ & $0.686$ & $0.686000$ & $0.00\%$ & $30.0$  & $30.000000$ & $0.00\%$ & $(0.0101,\ -0.1339)$ & satisfied \\
\bottomrule
\end{tabular}
}
\end{table}

\subsection{Computational Performance}
\label{sec:res_perf}

The factorization exists to remove the $\mathcal{O}(\dof^3)$ storage and arithmetic bottleneck of the dense Cartesian operator. The polyatomic construction enlarges the dense block of that operator from side $\kmax + 1$ to $N_{KI} = (\kmax+1)(\imax+1)$ (Section~\ref{sec:contraction}). We therefore count storage and measure contraction cost along both truncation axes, $\lmax$ and $\imax$, and test the ordering analysis of Section~\ref{sec:contraction} against measurement.

All measurements are taken at $\kmax = 2$ with the compiled contraction kernels of Section~\ref{sec:contraction}, executed on a single core of a $12$-core Apple M2 Pro ($32$~GB, macOS~14.4.1) under MATLAB R2023b. We report only relative performance: all four strategies run in the same environment, so the ratios below are independent of it. The timed quantity is one nonlinear operator evaluation $Q = \Ctens : (f \otimes f)$ including allocation of the output buffer, reported as the median over repeated calls, for four strategies: the dense Cartesian baseline, naive streaming, radial-first, and angular-first. All four kernels are nested loops over the geometric sparsity pattern and the dense radial--internal block, and none branches on a coefficient value. Their cost is therefore fixed by structure alone: the Gaunt nonzero count $\numG$, the channel count $\numT$, the angular size $\numQ$, and the block size $N_{KI}$. The timings accordingly use the true Gaunt geometry from the generator with the radial--internal blocks populated without running the nine-dimensional quadrature. This makes truncations reachable that a converged build could not afford, and no physical claim is made from these runs. The storage comparison is a count, taken over the same geometry:
\begin{equation}
    M_{\mathrm{dense}} = 8 \, \dof^3, \qquad
    M_{\mathrm{fact}} = 8 \left( \numT N_{KI}^3 + 5 \numG \right), \qquad
    \dof = \numQ N_{KI},
    \label{eq:storage_model}
\end{equation}
in bytes, with $\numQ = (\lmax+1)^2$. The factorized count carries the dense channel blocks of $\Rtens$ together with the coordinate-format routing list of $\Gtens$, at five doubles per transition (three index labels, the Gaunt weight, and the channel tag).

\subsubsection{Storage Compression}
\begin{figure}[t]
\centering
\includegraphics[width=\textwidth]{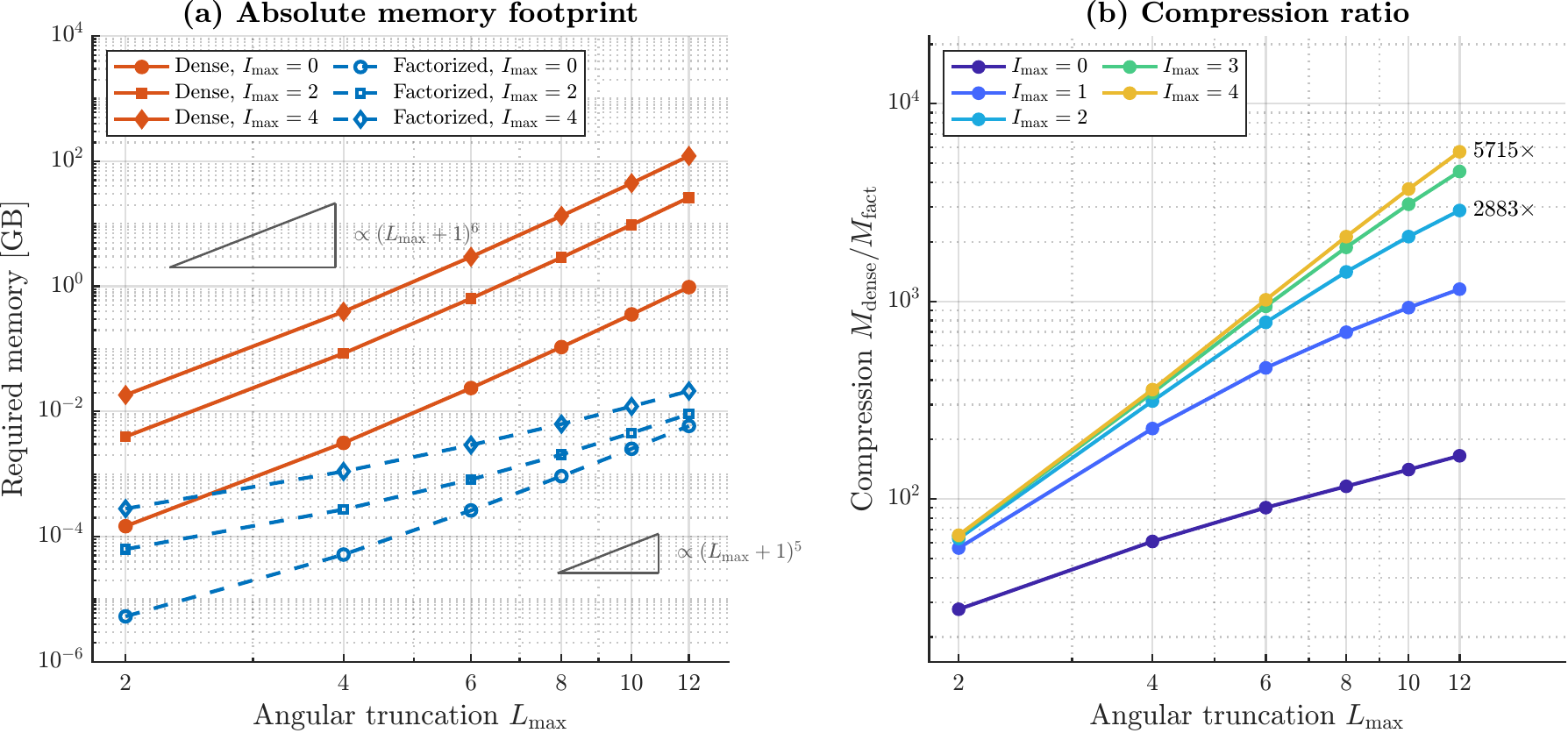}
\caption{Memory complexity of the factorized polyatomic collision operator at $\kmax = 2$, evaluated from the storage model \eqref{eq:storage_model}. (a)~Absolute footprint against the angular truncation for the dense Cartesian tensor and for the factorized representation at $\imax = 0$, $2$, and $4$. (b)~Compression $M_{\mathrm{dense}}/M_{\mathrm{fact}}$ against the angular truncation for $\imax = 0, \ldots, 4$, reaching $2883\times$ at $\lmax = 12$, $\imax = 2$ and $5715\times$ at $\imax = 4$.}
\label{fig:perf_storage}
\end{figure}

The two storage families separate over the whole measured range (Figure~\ref{fig:perf_storage}). The dense tensor follows $(\lmax+1)^6$ exactly, since $\dof^3 = (\lmax+1)^6 N_{KI}^3$, while the factorized footprint is governed by its two terms at different rates. The routing list grows as $(\lmax+1)^5$ and is independent of the internal truncation, whereas the channel blocks grow only as $(\lmax+1)^3$ and carry the whole internal dependence. Which term dominates is therefore set by $\imax$: at $\lmax = 12$ the routing list is $98\%$ of the factorized total for $\imax = 0$ but only $27\%$ for $\imax = 4$. Because the dense tensor grows faster than either term, the compression widens along both truncation axes at once.

At $\lmax = 12$ the compression ratio $M_{\mathrm{dense}}/M_{\mathrm{fact}}$ grows from $165\times$ at $\imax = 0$ to $5715\times$ at $\imax = 4$, climbing toward the purely geometric ceiling $\numQ^3/\numT = 7836$ as the dense blocks outgrow the routing list (Table~\ref{tab:perf_storage}).

\begin{table}[t]
\centering
\caption{Storage of the polyatomic collision operator at $\kmax = 2$: the dense Cartesian tensor of $\dof^3$ doubles against the factorized representation \eqref{eq:storage_model}, for angular truncations $\lmax = 2, \ldots, 12$ at the internal truncations $\imax = 2$ and $\imax = 4$. Here $\dof$ is the total number of degrees of freedom, $\numG$ the number of nonzero Gaunt transitions, and $\numT$ the number of physical channels. Storage is in binary units ($1$~GB $= 2^{30}$ bytes). The compression grows along both truncation axes and exceeds three orders of magnitude at the largest resolutions.}
\label{tab:perf_storage}
\begin{tabular}{@{}lcccccc@{}}
\toprule
$\lmax$ & $\dof$ & $\numG$ & $\numT$ & Dense [GB] & Factorized [MB] & Compression \\
\midrule
\multicolumn{7}{@{}l}{\emph{Internal truncation} $\imax = 2$} \\
$2$  & $81$   & $83$      & $11$  & $0.0040$  & $0.064$ & $63\times$ \\
$4$  & $225$  & $1{,}158$ & $42$  & $0.085$   & $0.278$ & $313\times$ \\
$6$  & $441$  & $6{,}460$ & $106$ & $0.639$   & $0.836$ & $783\times$ \\
$8$  & $729$  & $23{,}621$ & $215$ & $2.887$  & $2.097$ & $1410\times$ \\
$10$ & $1{,}089$ & $65{,}913$ & $381$ & $9.622$ & $4.633$ & $2127\times$ \\
$12$ & $1{,}521$ & $154{,}330$ & $616$ & $26.217$ & $9.313$ & $2883\times$ \\
\midrule
\multicolumn{7}{@{}l}{\emph{Internal truncation} $\imax = 4$} \\
$2$  & $135$  & $83$      & $11$  & $0.018$   & $0.286$ & $66\times$ \\
$4$  & $375$  & $1{,}158$ & $42$  & $0.393$   & $1.126$ & $357\times$ \\
$6$  & $735$  & $6{,}460$ & $106$ & $2.958$   & $2.976$ & $1018\times$ \\
$8$  & $1{,}215$ & $23{,}621$ & $215$ & $13.363$ & $6.437$ & $2126\times$ \\
$10$ & $1{,}815$ & $65{,}913$ & $381$ & $44.547$ & $12.325$ & $3701\times$ \\
$12$ & $2{,}535$ & $154{,}330$ & $616$ & $121.374$ & $21.749$ & $5715\times$ \\
\bottomrule
\end{tabular}
\end{table}

\subsubsection{Contraction Timings}
\begin{figure}[t]
\centering
\includegraphics[width=0.5\textwidth]{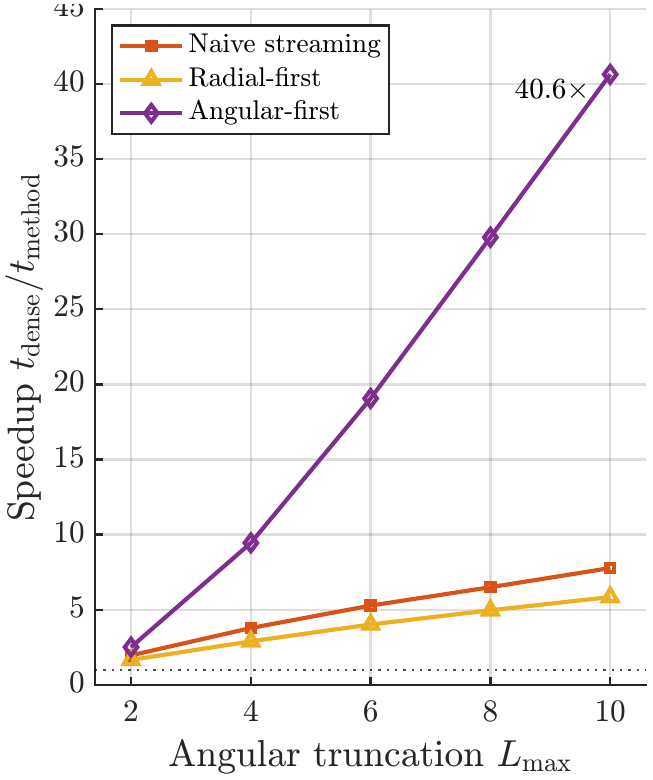}
\caption{Speedup of the three factorized contraction orderings over the dense Cartesian baseline at $\kmax = 2$, $\imax = 2$, against the angular truncation, timed single-core as one nonlinear operator evaluation including output allocation, median over repeated calls. The dense baseline is run up to $\lmax = 10$. The speedup grows monotonically with $\lmax$, with the angular-first ordering fastest at every truncation.}
\label{fig:perf_contraction}
\end{figure}

The contraction timings follow the same pattern (Figure~\ref{fig:perf_contraction}). At the largest dense-resident resolution, $\lmax = 10$ with $\imax = 2$, the angular-first ordering evaluates the operator $40.6\times$ faster than the dense Cartesian baseline, against $7.8\times$ for naive streaming and $5.8\times$ for the radial-first ordering. At $\lmax = 12$, where the baseline no longer fits, angular-first remains $7.6\times$ ahead of radial-first. The speedup tracks the $\mathcal{O}(\lmax^6)$ against $\mathcal{O}(\lmax^5)$ complexity gap and is consistent with the $37.2\times$ reported for the monatomic operator \citep{hiemstra2026wigner}.

The internal axis decides the ordering question of Section~\ref{sec:contraction} in favor of the angular-first ordering: the margin $t_{\mathrm{radial}}/t_{\mathrm{angular}}$ widens with internal resolution instead of narrowing. At $\lmax = 8$ it grows monotonically from $2.05$ to $8.16$ across $\imax = 0, \ldots, 4$, and at $\lmax = 6$ from $1.87$ to $6.35$ over the same range. Most of this margin is arithmetic rather than cache-related. The slice count falls from $\numS = 35$ at $\lmax = 2$ to $9100$ at $\lmax = 12$, i.e. from $0.42\numG$ to $0.06\numG$, so replacing $\numG$ by $\numS$ in the cubic term of Table~\ref{tab:contraction_orderings} predicts a margin of $4.9$ to $6.3$ over the measured range. The measured margins exceed that prediction by a further factor $1.23$ to $1.29$, which is the cache-locality contribution described by the analysis of Section~\ref{sec:contraction}. The angular-first intermediate stays first-level-cache resident over the whole range, at $1.8$~kB for $\imax = 4$, whereas the radial-first intermediate of $\numG N_{KI}$ entries grows with the internal resolution ($2.7$~MB at $\lmax = 8$, $\imax = 4$). The enlarged radial--internal block therefore penalizes the radial-first ordering alone.

\subsubsection{Feasibility Boundaries}
The dense baseline is not run at $(\lmax, \imax) = (12, 2)$, $(8, 4)$, and $(10, 4)$, where its footprint \eqref{eq:storage_model} reaches $26.2$, $13.4$, and $44.5$~GB. The factorized operator evaluates all three in under $13$~MB. A dense implementation carrying $32$-bit linear indices would stop earlier still, at $\dof > 1290$, where $\dof^3$ overflows the signed $32$-bit range. Our kernels index in $64$-bit, so memory, not indexing, is the binding constraint.

\section{Conclusion}

We presented a conservative spectral method for the Boltzmann collision operator of polyatomic gases in the continuous-internal-energy framework, extending the Wigner--Eckart factorization of the monatomic companion \citep{hiemstra2026wigner} to molecular gases. Because the internal energies and partition parameters are scalars under spatial rotations, the $\mathrm{SO}(3)$ reduction survives the internal coupling intact: the twelve-dimensional collision integral collapses onto a nine-dimensional kinematic core routed by the unchanged sparse geometric Gaunt tensor.

The implementation is defined by five methodological contributions:
\begin{itemize}
  \item \textbf{Dimensional Reduction:} Extending the geometric factorization to the Borgnakke--Larsen internal-energy exchange, reducing the twelve-dimensional polyatomic collision integral to a nine-dimensional physical core with the monatomic angular geometry reused verbatim.
  \item \textbf{Spectral Quadrature:} Composing exactness-matched Gauss rules with an auxiliary Laplace representation of the fractional energy couplings on both collision channels, so that every integration direction converges spectrally: to machine plateaus below $10^{-13}$ at integer interaction exponents, and to below $2 \times 10^{-11}$ at the fractional exponents of the production gases.
  \item \textbf{Invariant Conservation:} Extending the null-space embedding to the polyatomic invariants, where an orthogonal projection removes the two-mode energy invariant while preserving the translational--internal exchange: mass is conserved exactly and total energy to below $3 \times 10^{-14}$ under nonlinear evolution.
  \item \textbf{Storage Compression:} Growing the compression of the factorized operator along both truncation axes, from three to nearly four orders of magnitude ($2883\times$ to $5715\times$) over the dense Cartesian tensor at the largest resolutions measured.
  \item \textbf{Execution Acceleration:} Carrying the cache-optimized angular-first contraction to the enlarged radial--internal blocks, evaluating the operator $40.6\times$ faster than the dense baseline, with a margin that widens as the internal resolution grows.
\end{itemize}

The framework was validated against the monatomic operator, reproduced exactly at the internal collapse point and approached continuously in the internal dimension. It was validated further against the Landau--Teller relaxation law and the Wang Chang--Uhlenbeck spectrum for Maxwell molecules, against the closed-form frozen Prandtl number of \citet{djordjic:2023}, matched to $2 \times 10^{-7}\%$, and against the classical hard-sphere viscosity factor $205/202$. These comparisons are parameter-free. At the published kernel parameters the operator reproduces the transport coefficients of N$_2$, CO, and H$_2$ with nothing fitted, but the agreement with experiment there is inherited from the calibration of \citet{djordjic:2023}, whose parameters were fitted to those same two coefficients through the same closure. That comparison is a code-to-code verification: it establishes that a nine-dimensional quadrature reproduces their symbolically evaluated production coefficients to within $0.0002\%$. Calibrating the two free parameters within our operator recovers the published pair as well, so the symbolic evaluation and the quadrature agree on the parameters as well as on the coefficients. The study also showed that the bulk-to-shear ratio, unlike the Prandtl number, is sensitive to the order at which the transport coefficients are extracted. By removing the $\mathcal{O}(\kmax^3 \imax^3 \lmax^6)$ storage and arithmetic bottlenecks of dense polyatomic formulations, the factorization provides a foundation for high-resolution deterministic kinetic simulation of molecular gases.

\appendix
\section{Sum-Factorized Offline Assembly of the Physical Tensor}
\label{app:assembly}

This appendix organizes the offline build: for each angular channel $\chan = (l_1, l_2, l_3)$ admitted by the Wigner selection rules, the quadrature directions of Section~\ref{sec:quadfamilies} combine into the dense block $\Rtens_{k_1 k_2 k_3, i_1 i_2 i_3}^{l_1 l_2 l_3}$, at a cost incurred once per kernel calibration. The online contraction (Section~\ref{sec:contraction}) is untouched. Factoring the polar deflection outside the inner integrations defines the inner polyatomic block
\begin{equation}
    \polyblock^{(\chan)}(\cmenergy, \fracrel, \halfangle) = \int_{-1}^{1} \Big[ \polyblock_{\mathrm{gain}}^{(\chan)}(\cmenergy, \fracrel, \halfangle, \cos\defangle) - \polyblock_{\mathrm{loss}}^{(\chan)}(\cmenergy, \fracrel, \halfangle, \cos\defangle) \Big] \, d(\cos\defangle),
    \label{eq:inner_block}
\end{equation}
where each one-dimensional integral below is evaluated by the corresponding tabulated rule, whose mapped weight we denote generically by $W$. The spectral indices carried by each stage are suppressed for clarity.

Because parts of the kernel are already carried by the quadrature weights, the integrand of every stage evaluates the \emph{reduced} kernel $\redker$, obtained from \eqref{eq:dsmc_kernel} by stripping the weight-folded factors: the energy power $\cmenergy^{\vhs/2}$ (energy weight; $\cmenergy^{\vhs}$ on the frozen-spectral pass below), the cone power $\fracrel^{\vhs}$ (respectively $\halfangle^{\vhs}$) absorbed by the radial weight $\radweight^{(\chan)}$ below, and the partition power $\partkin^{\vhs/2}$ (kinetic-partition weight). On the modulated branches the bounded partition prefactors $\partint^{\modexp}$, $(1-\partint)^{\modexp}$, and $(1-\partkin)^{\modexp}$ are likewise weight-folded, into the per-channel shifted Gauss--Jacobi rules of Section~\ref{sec:laplace}. Of the factor $\urel^{\vhs}$, only the bounded per-patch remainder (e.g., $2^{\vhs/2} (1 + (1-\fracrel^2)\duffy^2)^{\vhs/2}$ on Patch 1, cf.~\eqref{eq:duffy_speeds}) survives pointwise in $\redker$.

The loss integrand evaluates only pre-collision states. It is therefore independent of the azimuthal angle and of the geometric reconstruction \eqref{eq:post_reconstruction}, so the azimuthal integral evaluates to $2\pi$ and the loss filter collapses to its zonal form (Section~\ref{sec:reduction}):
\begin{equation}
    \polyblock_{\mathrm{loss}}^{(\chan)}(\cmenergy, \fracrel, \halfangle, \cos\defangle) = 2\pi \, P_{\mathrm{loss}}^{(\chan)}(\halfangle) \, \radbasis_{k_1, l_1}(\vmag) \int_0^\infty \!\! \int_0^\infty \!\! \int_0^1 \!\! \int_0^1 W_{\mathrm{4D}} \, \redker \, \intbasis_{i_1}(\Ipre) \intbasis_{i_2}(\Ipre) \intbasis_{i_3}(\Jpre) \, d\partint \, d\partkin \, d\Jpre \, d\Ipre,
    \label{eq:loss_block}
\end{equation}
where $W_{\mathrm{4D}}$ collects the four internal weight functions of the quadrature table and the four inner integrals tensorize into independent one-dimensional sums. The deflection angle enters \eqref{eq:loss_block} only through the kernel. For the isotropic kernel \eqref{eq:dsmc_kernel} the loss block is constant in $\cos\defangle$.

The gain integrand requires the post-collision states $\vmag'$, $\vpostunit$, and $\Ipost$. Applying Fubini's theorem along the dependency chain of Section~\ref{sec:kinematics} nests the integrals so that the three-dimensional scattering geometry is never re-evaluated inside the internal-energy loops:
\begin{equation}
\begin{aligned}
    \stageM(\cmenergy, \fracrel, \halfangle, \cos\defangle, \Ipre, \Jpre, \partkin) &= \int_0^1 W(\partint) \, \redker \, \intbasis_{i_1}\big(\Ipost\big) \, d\partint, \\
    \stageS^{(\chan)}(\cmenergy, \fracrel, \halfangle, \cos\defangle, \Ipre, \Jpre, \partkin) &= \int_0^{2\pi} \radbasis_{k_1, l_1}(\vmag') \, P_{\mathrm{gain}}^{(\chan)}(\vpostunit_0, \vpreunit_0, \wpreunit_0) \, d\aziangle, \\
    \stageN^{(\chan)}(\cmenergy, \fracrel, \halfangle, \cos\defangle, \Ipre, \Jpre) &= \int_0^1 W(\partkin) \, \stageM \, \stageS^{(\chan)} \, d\partkin, \\
    \stageO^{(\chan)}(\cmenergy, \fracrel, \halfangle, \cos\defangle, \Ipre) &= \int_0^\infty W(\Jpre) \, \intbasis_{i_3}(\Jpre) \, \stageN^{(\chan)} \, d\Jpre, \\
    \polyblock_{\mathrm{gain}}^{(\chan)}(\cmenergy, \fracrel, \halfangle, \cos\defangle) &= \int_0^\infty W(\Ipre) \, \intbasis_{i_2}(\Ipre) \, \stageO^{(\chan)} \, d\Ipre.
\end{aligned}
\label{eq:gain_nesting}
\end{equation}
The innermost stage $\stageM$ carries only the kernel and the post-collision internal basis, because $\vmag'$ and $\vpostunit$ are independent of $\partint$. The azimuthal stage $\stageS^{(\chan)}$ encapsulates the geometric reconstruction and the $\mathcal{O}(\kmax)$ radial polynomial evaluation. The two are merged over their shared dependency $\partkin$ and subsequently contracted against the incident and target internal bases. This ordering also resolves the square root $\urel' = \sqrt{2 \partkin \Etot}$ analytically. The trapezoidal average over $\aziangle$ annihilates the $\aziangle$-odd components of $\scatterdir$, and the symmetry of the deflection rule in $\cos\defangle$ cancels the remaining odd powers, so only even powers $(\urel')^{2n} = (2 \partkin \Etot)^n$ contribute to the assembled block. The surviving integrands are polynomial in $\Ipre$, $\Jpre$, and $\partkin$, the premise underlying the internal exactness bounds of Section~\ref{sec:quadfamilies}. On the modulated and frozen branches the same holds against the shifted and rescaled weights of Section~\ref{sec:laplace}, whose auxiliary representation restores the product structure.

The inner block \eqref{eq:inner_block} depends only on the spatial triple $(\cmenergy, \fracrel, \halfangle)$ and integrates into the two Duffy patches of the monatomic spatial architecture:
\begin{equation}
\begin{aligned}
    \Rtens_{1}^{(\chan)} &= \int_0^\infty \!\! d\cmenergy \, \cmenergy^{\vhs/2} e^{-\cmenergy} \int_0^1 \!\! d\fracrel \, \radweight^{(\chan)}(\cmenergy, \fracrel) \, \fracrel \int_0^1 \!\! d\duffy \, \polyblock^{(\chan)}(\cmenergy, \fracrel, \fracrel\duffy), \\
    \Rtens_{2}^{(\chan)} &= \int_0^\infty \!\! d\cmenergy \, \cmenergy^{\vhs/2} e^{-\cmenergy} \int_0^1 \!\! d\halfangle \, \halfangle \int_0^1 \!\! d\duffy \, \polyblock^{(\chan)}(\cmenergy, \halfangle\duffy, \halfangle) \, \radweight^{(\chan)}(\cmenergy, \halfangle\duffy),
\end{aligned}
\label{eq:patch_assembly}
\end{equation}
where the radial weight $\radweight^{(\chan)}$ collects the pre-collision radial polynomials $\radbasis_{k_2, l_2}(\vmag) \radbasis_{k_3, l_3}(\wmag)$, the cone power of the kernel, and the mapped spatial measure. Its closed form coincides, up to the polyatomic reference normalization, with Appendix~A, Eqs.~(47)--(49) of \citet{hiemstra2026wigner}. On Patch 1 ($\halfangle = \fracrel \duffy$) the radial weight depends only on the outer coordinate and factors out of the inner Duffy loop; on Patch 2 ($\fracrel = \halfangle \duffy$) it depends on the inner coordinate and is evaluated inside. The nesting \eqref{eq:gain_nesting} bounds the arithmetic: the $\mathcal{O}(\kmax)$ radial evaluation resides in $\stageS^{(\chan)}$ and is factored out of the innermost $\partint$ loop, so each stage touches every quadrature node exactly once per retained spectral index it carries.

The frozen channel of \eqref{eq:dsmc_kernel} is assembled in a separate elastic pass. Its Dirac factors are normalized against the weighted partition measure,
\begin{equation}
    \int_0^1 \!\! \int_0^1 \delta(\partkin - \bar{\partkin}) \, \delta(\partint - \bar{\partint}) \, \phipart(\partint) \psipart(\partkin) (1-\partkin) \partkin^{1/2} \, d\partint \, d\partkin := 1,
    \label{eq:frozen_normalization}
\end{equation}
which consumes the partition functions and the collision Jacobian of the measure, but not the density $\partkin^{\vhs/2}$ carried by the kernel itself. The frozen channel therefore runs without a partition measure, as the elastic process of \eqref{eq:dsmc_kernel} weighted by the pinned partition factor $\bar{\partkin}^{\vhs/2} = (\urel^2/2\Etot)^{\vhs/2}$. The build runs the monatomic spatial architecture with $\urel' = \urel$ and the internal energies transported unchanged, producing in place of the four-dimensional inner integrals of \eqref{eq:gain_nesting} the internal overlap integrals $\int_0^\infty W(\Ipre) \intbasis_{i_1}(\Ipre) \intbasis_{i_2}(\Ipre) \, \redker \, d\Ipre$ and $\int_0^\infty W(\Jpre) \intbasis_{i_3}(\Jpre) \, \redker \, d\Jpre$. The pinned partition factor is not evaluated pointwise. Following Section~\ref{sec:laplace}, its split $\bar{\partkin}^{\vhs/2} = (\urel^2/2)^{\vhs/2} \, \Etot^{-\vhs/2}$ folds the translational factor into the kernel power (the frozen-spectral integrand carries $\urel^{2\vhs}$, desingularized against the rekeyed energy weight $\cmenergy^{\vhs} e^{-\cmenergy}$), while the coupled factor $\Etot^{-\vhs/2}$ passes through the auxiliary Laplace integral. On the modulated branch it combines with the power $\Etot^{-\modexpf}$ into the single exponent $p = \vhs/2 + \modexpf$: one auxiliary integral, not two. Under the auxiliary rule the overlap integrands are again polynomial against their own Laguerre weights, so the internal quadrature is exact. For the base frozen kernel at $\vhs = 0$ this factor is unity, no auxiliary integral is required, and orthonormality reduces the overlaps to the Kronecker structure $\delta_{i_1 i_2}$, $\delta_{i_3 0}$.

The auxiliary correction of Section~\ref{sec:laplace} slots into this pipeline as follows. At each auxiliary node $\auxmap_q$, the internal-energy nodes are rescaled by $\auxscale = 1 - \auxmap_q$ with the active axis carrying the shifted Laguerre parameter $\polyalpha + \modexp$. The partition prefactors are carried by the per-channel shifted Gauss--Jacobi rules of Section~\ref{sec:laplace}, with the $\Ipre$-active and $\Jpre$-active channels summed. The Gaussian factor $e^{-\auxvar \urel^2/2}$ multiplies the reduced kernel $\redker$ pointwise at the spatial nodes (with the factors held by the shifted weights likewise stripped from $\redker$). The builds are combined with the Gauss--Jacobi weights of \eqref{eq:aux_jacobi}, which carry the node-scaling powers stripped from the builds. Each channel then assembles as base plus amplitude-scaled correction, $\Rtens_{\mathrm{base}} + \modamp \, \Rtens_{\mathrm{corr}}$ (with $\modampf$ on the frozen channel), and the frozen channel's auxiliary passes run through the same slot with the rekeyed energy weight and the combined exponent $p$. The implementation makes two auxiliary calls per node, one $\Ipre$-active and one $\Jpre$-active, so the offline cost of the extended kernel is $1 + 2\numaux$ base builds per kernel calibration: the nine-dimensional non-frozen build once, plus two rescaled evaluations at each of the $\numaux$ auxiliary nodes. The frozen channel's auxiliary passes run on the cheaper elastic architecture.

The storage footprint mirrors the monatomic analysis: the dense physical tensor holds $(\kmax+1)^3 (\imax+1)^3$ entries for each of the $\numT \approx \tfrac{1}{4}\lmax^3$ angular channels admitted by the triangle and parity selection rules (Eq.~(35) of \citet{hiemstra2026wigner}), and the sparse Gaunt routing adds $\mathcal{O}(\lmax^5)$ nonzero geometric transitions. The offline build evaluates the nested stages of \eqref{eq:gain_nesting} once per retained index combination over the tensorized grid of Table~\ref{tab:quadrature_bounds}.

\section{Analytic Relaxation Rates for a Maxwell Linear Molecule}
\label{app:maxwell}

This appendix derives the exact relaxation eigenvalues that anchor the Maxwell-kernel validations of Section~\ref{sec:res_relax}, including the Landau--Teller rate \eqref{eq:landau_teller} in terms of the internal degrees of freedom. We consider the base Borgnakke--Larsen Maxwell kernel ($\vhs = 0$, $\inelprob = 1$, $\modamp = \modampf = 0$) for a molecule with $\dofpoly$ internal degrees of freedom ($\polyalpha = \dofpoly/2 - 1$) alongside $d_t = 3$ translational degrees of freedom. A colliding pair redistributes its energy over $d_t + 2\dofpoly$ active degrees of freedom in the center-of-mass frame. The linear molecule of the validation runs, $\dofpoly = 2$ ($\polyalpha = 0$), gives $7$ such degrees of freedom.

\subsection{Shear Relaxation Rate}
The shear rate is governed by the deviatoric stress mode $\phi_{\mathrm{shear}} \propto \{ \mathbf{v} \otimes \mathbf{v} \}$, the $(k, l) = (0, 2)$ mode of Section~\ref{sec:res_mono}. Because this mode carries no internal energy, it is unaffected by the partitions $\partkin$ and $\partint$. Integrating the post-collision gain term over the isotropic scattering solid angle yields an expected gain fraction of $1/2$; subtracting the unit loss term,
\begin{equation}
    \lambda_{\mathrm{shear}} \propto \left( \frac{1}{2} - 1 \right) = -\frac{1}{2}.
    \label{eq:lambda_shear}
\end{equation}
In units of the equilibrium collision frequency this is $\lambda_{\mathrm{shear}} = -\nu_0/2$, independently of $\dofpoly$. Section~\ref{sec:res_mono} verifies this value against the discrete operator, together with the decoupling from the internal sector on which it rests.

\subsection{Energy-Exchange Relaxation Rate}
The bulk viscosity is dictated by the thermal relaxation between the translational and internal energy pools, carried by the isotropic mode $\phi_{\mathrm{exch}} \propto \Ipre - \tfrac{\dofpoly}{d_t} E_{\mathrm{kin}}$ orthogonal to the total-energy invariant (the exchange mode preserved by the conservation projection \eqref{eq:energy_projection}). In an inelastic collision the expected fraction of the total energy deposited into the internal modes follows from the Beta distribution of the kinetic partition $\partkin$ under the weight of Table~\ref{tab:quadrature_bounds}:
\begin{equation}
    \langle 1 - \partkin \rangle = \frac{2\dofpoly}{d_t + 2\dofpoly}.
\end{equation}
Projecting the gain operator onto $\phi_{\mathrm{exch}}$ with this expectation yields the prefactor $\dofpoly / (d_t + 2\dofpoly)$; subtracting the loss term,
\begin{equation}
    \lambda_{\mathrm{exch}} \propto \left( \frac{\dofpoly}{d_t + 2\dofpoly} - 1 \right) = -\frac{d_t + \dofpoly}{d_t + 2\dofpoly}.
    \label{eq:lambda_bulk}
\end{equation}
In units of the equilibrium collision frequency, and with $d_t = 3$, this is the Landau--Teller rate $\lambda_{\mathrm{LT}} = (3 + \dofpoly)/(3 + 2\dofpoly)$ used in \eqref{eq:landau_teller}. It is a property of the Borgnakke--Larsen energy partition rather than of the exponential relaxation law itself, and it takes the value $5/7$ at $\dofpoly = 2$.

\subsection{Eigenvalue Ratio}
In the ratio of the exchange to shear relaxation eigenvalues the dimensional collision-frequency constants cancel, giving at $\dofpoly = 2$
\begin{equation}
    \frac{\lambda_{\mathrm{exch}}}{\lambda_{\mathrm{shear}}} = \frac{-5/7}{-1/2} = \frac{10}{7}.
    \label{eq:ratio_bulk_shear}
\end{equation}
This ratio is basis- and normalization-independent. It is the analytic benchmark for the energy-exchange eigenvalue reported in the text of Section~\ref{sec:res_mono}.



\section*{Acknowledgements}
R.\,R. Hiemstra was funded by the European Union (Marie Skłodowska-Curie grant agreement No. 101105786). Views and opinions expressed are however those of the author(s) only and do not necessarily reflect those of the European Union or the European Research Executive Agency. Neither the European Union nor the granting authority can be held responsible for them.

\section*{Computer Code Availability}
The official software codebase is openly available on \href{https://github.com/simkinetic/collision-factorization/tree/v2.0.0}{GitHub}. The specific version used to generate the results in this paper is permanently archived on Zenodo under the DOI: \href{https://doi.org/10.5281/zenodo.22177473}{10.5281/zenodo.22177473}~\citep{kessler_zenodo_2026}.

\printcredits

\bibliographystyle{cas-model2-names}

\bibliography{cas-refs}

\end{document}